\documentclass[12pt]{article}%
\usepackage{amsmath}
\usepackage{amsfonts}
\usepackage{amssymb}
\usepackage{graphicx}
\usepackage{subcaption}
\usepackage{xcolor}
\usepackage{zymacros}%
\usepackage[ruled, vlined]{algorithm2e}

\begin{document}

\title{Linearized Learning Methods with Multiscale Deep Neural Networks for Stationary Navier-Stokes Equations with Oscillatory Solutions}
\author{Lizuo Liu\thanks{Department of Mathematics, Southern Methodist
University, Dallas, TX 75275.} \and Bo Wang \thanks{LCSM(MOE), School of Mathematics and Statistics, Hunan Normal University, Changsha, Hunan, 410081, PR China.}
\and Wei Cai\thanks{Corresponding author, Department of Mathematics, Southern
Methodist University, Dallas, TX 75275(\texttt{cai@smu.edu}).} }
\maketitle

\begin{abstract}
   In this paper, we present linearized learning methods to accelerate the convergence of training for stationary nonlinear Navier-Stokes equations. To solve the stationary nonlinear Navier-Stokes (NS) equation, we integrate the procedure of linearization of the nonlinear convection term in the NS equation into the training process of multi-scale deep neural network approximation of the NS solution. Four forms of linearizations are considered. After a benchmark problem, we solve the highly oscillating stationary flows utilizing the proposed linearized learning with multi-scale neural network for complex domains. The results show that multiscale deep neural network combining with the linearized schemes can be trained fast and accurately.
\end{abstract}

\section{Introduction}
\indent
Deep neural network (DNN) machine learning has been studied as alternative numerical methods for solving partial differential equations arising from fluid dynamics \cite{Wang_2020} and wave propagations \cite{phasednn}, which has the potential of a flexible meshless method to solve governing equations from fluid and solid mechanics applications in complex geometries as an alternative method to traditional finite element method. Moreover, it has shown much power in handling high dimensional parabolic PDEs \cite{raissi,han, zhang}. Recent work on DNNs have shown that they have a frequency dependence performance in learning solution of PDEs and fitting of functions. Namely, the lower frequency components of the solution are learned first and quickly compared with the higher frequency components \cite{Xu_2020}. Several attempts have been made to remove such a frequency bias for the DNN performances. One of the ideas is to convert the higher frequency content of the solution to a lower frequency range so the conventional DNN can learn the solution quickly. One way to achieve this goal is to use a phase shift\cite{phasednn} while the other way is to introduce a multiscale structure of the DNN \cite{CiCP-28-1970} where neural networks with different scales will target different ranges of the frequency in the solutions. The PhaseDNN has been shown to be very effect for high frequency wave propagations while the MscaleDNN has been used to learn highly oscillatory Stokes flow solutions in complex domains as well as high dimensional PDEs \cite{zhang}.

Additional difficulties arise when if there are nonlinearity in the PDEs. Based on the results from Jin et al.\cite{nsfnet}, it is known that it might cost \(\text{O}\left( 10^{4} \right)\) epochs to solve a benchmark problem, which is ineffective and impractical especially when highly oscillating problems are to be considered. Also, it is found that the MscaleDNN applied to the nonlinear Navier-Stokes equation did not produce the dramatic improvement over convential DNNs as in the case of the linear Stoke equations\cite{Wang_2020}.

The learning of the solution of linear PDEs via least square residual of the PDEs is in some sense equivalent to a fitting problems in the frequency domain from the Parsevel equality of the Fourier transform. So it is understandable to see a multiscale DNN’s performance improvement for a fitting problem also holds for the solution of linear PDEs. Considering the MscaleDNN’s good performance for the linear Stokes problem, we will develop in this paper a linearized learning procedure for the Navier-Stokes equations by integrating a linearization of the Navier-Stokes equation in the loss function of the MscaleDNN and dynamically updating the linearization as the learning is being carried out. Numerical results have demonstrated the convergence of this approach, which has produced highly accurate approximation to oscillatory solutions of the Navier-Stokes equations.

The rest of the paper will be organized as follows. Section \ref{sec:Multiscale} will introduce the multiscale DNN structure and then four linearized learning for the Navier-Stokes equation will be proposed in Section \ref{sec:linear}. Numerical tests of the linearized learning schemes will be conducted for a 2-D oscillatory flows in a domain containing one or multiple random cylinder(s) in Section \ref{sec:results}. Finally, conclusion and future work will be discussed in Section \ref{sec:future}.

\noindent

\section{Linearized iterative method for stationary Navier–Stokes equation}
Consider the following stationary Navier-Stokes equation
\begin{equation}
	\left\{\begin{array}{ll}
		(\boldsymbol{u} \cdot \nabla) \boldsymbol{u}- \nu \Delta \boldsymbol{u} + \nabla p =\boldsymbol{f} & \text { in } \Omega,  \\
		\nabla \cdot \boldsymbol{u}=0 & \text { in } \Omega, \\
		\boldsymbol{u}=\boldsymbol{g} & \text { on } \Gamma_{D}, \\
	\end{array}\right.
	\label{NS}
\end{equation}
where $\Omega$ is an open bounded domain in \(\mathbb{R}^d, d=2, 3\).

To solve the Stationary Navier-Stokes equation in a least square approach, intuitively, the PDEs (\ref{NS}) are written in a system of first order equations as in \cite{Wang_2020}. A velocity-gradient of velocity-pressure formulation will be used for the Navier Stokes equations.

\begin{subequations}
    \begin{align}
            -\nu \nabla \cdot \underline{\mathbf{U}}+\underline{\mathbf{U}} \cdot \mathbf{u}+\nabla p &=\mathbf{f} & & \text { in } \Omega \label{eq_1}\\
            \underline{\mathbf{U}}-(\nabla \mathbf{u})^{\top} &=0 & & \text { in } \Omega \label{eq_2}\\
            \nabla \cdot \mathbf{u} &=0 & & \text { in } \Omega \label{eq_3}
     \end{align}
\label{lsNS}
\end{subequations}

To obtain the governing equation for the pressure, the Navier Stokes equation is reformulated to Poisson equation with respect to pressure by taking divergence on both sides of equation \(\left( \ref{eq_1} \right)\) and applying equation \(\left( \ref{eq_3} \right)\),
\begin{equation}
         \Delta p + 2(-u_{x}v_{y}+ u_{y}v_{x})= \nabla\cdot \mathbf{f}.
\label{PressurePossion}
\end{equation}
Then, the loss function can be defined as
\begin{equation}
    \begin{aligned}
				L_{V g V P}\left(\theta_{u}, \theta_{p}, \theta_{\underline{\mathbf{U}}} \right):=&\|\nu \nabla \cdot {\underline{\mathbf{U}}} - \underline{\mathbf{U}}\cdot \boldsymbol{u}-\nabla p+\boldsymbol{f}\|_{\Omega}^{2}\\
        &+\alpha\|\Delta p+2(-u_{x}v_{y}+u_{y}v_{x}) -\nabla \cdot \boldsymbol{f}\|_{\Omega}^{2}\\
				&+\|\nabla \boldsymbol{u}-\underline{\mathbf{U}}\|_{\Omega}^{2} \\
        &+\|\nabla \cdot \boldsymbol{u}\|_{\Omega}^{2}\\
        &+\beta\|\boldsymbol{u}-\boldsymbol{g}\|_{\partial \Omega}^{2}.
    \end{aligned}
    \label{loss_nl}
\end{equation}
where \(\alpha\) and \(\beta\) are the penalty terms.
On the contrary, unlike the results in the paper \cite{Wang_2020}, it will be shown that the training of the network based on the formulation  \(\left( \ref{lsNS} \right)\) including the Nonlinear first order system and  the Poisson equation  \(\left( \ref{PressurePossion} \right)\) converges slowly even applying the MscaleDnn, as shown in section \ref{results}.

In the paper \cite{HE20091351}, three iterative methods were studied for solving the stationary Navier-Stokes equations. Assume \(u_h^n\) be the solution at n-th step of the iteration, three iterative schemes are given as follows.
\begin{itemize}
	\item \textbf{Iterative method I}
	\begin{equation}
		a\left(u_{h}^{n}, v_{h}\right)-d\left(v_{h}, p_{h}^{n}\right)+d\left(u_{h}^{n}, q_{h}\right)+b\left(u_{h}^{n-1}, u_{h}^{n-1}, v_{h}\right)=\left(f, v_{h}\right),
	\end{equation}
	\item \textbf{Iterative method II}
	\begin{equation}
		\begin{array}{l}
			a\left(u_{h}^{n}, v_{h}\right)-d\left(v_{h}, p_{h}^{n}\right)+d\left(u_{h}^{n}, q_{h}\right)
			+b\left(u_{h}^{n}, u_{h}^{n-1}, v_{h}\right)\\+b\left(u_{h}^{n-1}, u_{h}^{n}, v_{h}\right)
			=b\left(u_{h}^{n-1}, u_{h}^{n-1}, v_{h}\right)+\left(f, v_{h}\right),
		\end{array}
	\end{equation}
	\item \textbf{Iterative method III}
	\begin{equation}
		a\left(u_{h}^{n}, v_{h}\right)-d\left(v_{h}, p_{h}^{n}\right)+d\left(u_{h}^{n}, q_{h}\right)+b\left(u_{h}^{n-1}, u_{h}^{n}, v_{h}\right)=\left(f, v_{h}\right),
	\end{equation}
\end{itemize}
where \(a\left( u,v \right)=\nu\left( \nabla u , \nabla v \right),u,v \in X\),
\(d\left( v,q \right)=\left( q,\text{div} v \right), v\in X, q \in M\),
\(b\left( u,v,w \right) = \left( \left( u\cdot \nabla  \right)v + \frac{1}{2} \text{uv},w \right), u,v,w\in X\). The \(X, Y, M\) are Hilbert spaces \(X=H_{0}^{1}(\Omega)^{2},  Y=L^{2}(\Omega)^{2}\),\(M=L_{0}^{2}(\Omega)\).

Assume \(a_1(u,v,w) = \left( \left( u\cdot \nabla \right)v,w \right)\), then \(b\left( u,v,w \right) = \frac{1}{2}a_1\left( u,v,w \right) - \frac{1}{2}a_1\left( u,w,v \right)\). Obviously, the bilinear term \(a\left( \cdot,\cdot \right)\) is continuous and coercive on \(X\times X\); the bilinear \(d\left( \cdot, \cdot \right)\) satisfies that for all \(q\in M\)
\begin{equation}
	sup_{v\in X} \frac{\|d\left( v,q \right)\|}{\left\| \nabla v\right\|_{0}} \geq \beta_0 \|q\|_0,
\end{equation}
where \(\beta_{0}>0\). The trilinear form \(a_1\left( \cdot, \cdot, \cdot \right)\) satisfies
\begin{equation}
	|a_1\left( u,v,w \right)| \leq N \|\nabla u\|_{0} \|\nabla v\|_{0} \|\nabla w\|_{0},
	\label{eq:a_1}
\end{equation}
where \(N>0\).

Under the stability condition
\begin{equation}
	\frac{4N \|f\|_{-1}}{\nu^2}<1, \frac{25N \|f\|_{-1}}{3\nu^2}<1
	\label{stability}
\end{equation}
and the uniqueness condition
\begin{equation}
	\frac{N \|f\|_{-1}}{\nu^2}<1,
	\label{uniqueness}
\end{equation}
the following error estimates for all three schemes can be obtained.

Given the mesh size \(h\) and iterative step \(m\), for iterative methods I and III, we have
\begin{equation}
	\begin{aligned}
		&\nu\left\|u-u_{h}^{m}\right\|_{0} \leqslant C h^{2}+C v\left\|u_{h}^{m}-u_{h}^{m-1}\right\|_{0}, \\
		&\nu\left\|\nabla\left(u-u_{h}^{m}\right)\right\|_{0}+\left\|p-p_{h}^{m}\right\|_{0} \leqslant C h+C v\left\|u_{h}^{m}-u_{h}^{m-1}\right\|_{0}.
	\end{aligned}
\end{equation}

For iterative method II, we have
\begin{equation}
	\begin{aligned}
		&\nu\left\|u-u_{h}^{m}\right\|_{0} \leqslant C h^{2}+c|\log h|^{1 / 2}\left\|\nabla\left(u_{h}^{m}-u_{h}^{m-1}\right)\right\|_{0}\left\|u_{h}^{m}-u_{h}^{m-1}\right\|_{0}, \\
		&\nu\left\|\nabla\left(u-u_{h}^{m}\right)\right\|_{0}+\left\|p-p_{h}^{m}\right\|_{0} \leqslant C h+\\
		&\quad\quad\quad\quad\quad\quad\quad c|\log h|^{1 / 2}\left\|\nabla\left(u_{h}^{m}-u_{h}^{m-1}\right)\right\|_{0}\left\|u_{h}^{m}-u_{h}^{m-1}\right\|_{0}.
	\end{aligned}
\end{equation}

The strong form of these three iterative methods are given below.
\begin{itemize}
	\item {\bf Iterative method I}
	\begin{subequations}
		\begin{align}
			-\nu \Delta u^{n}  + u^{n-1}_{x}u^{n-1}+u^{n-1}_{y}v^{n-1} + p_{x} &=\mathbf{f}_{1},& & \\
			-\nu \Delta v^{n}  +  v^{n-1}_{x}u^{n-1}+v^{n-1}_{y}v^{n-1} + p_{y} &=\mathbf{f}_{2}.& &
		\end{align}
		\label{eq:str_scheme1}
	\end{subequations}
	\item {\bf Iterative method II}
	
	\begin{subequations}
		\begin{align}
			\begin{split}
				-\nu \Delta u^{n}  + (u^{n}_{x}u^{n-1}+u^{n}_{y}v^{n-1} &+ u^{n-1}_{x}u^{n}+u^{n-1}_{y}v^{n} )+ p_{x} \\
				&=\mathbf{f}_{1} + u^{n-1}_{x}u^{n-1}+u^{n-1}_{y}v^{n-1},
			\end{split}
		\end{align}
		\begin{align}
			\begin{split}
				-\nu \Delta v^{n}  + ( v^{n}_{x}u^{n-1}+v^{n}_{y}v^{n-1}&+v^{n-1}_{x}u^{n}+v^{n-1}_{y}v^{n} )+ p_{y} \\
				&=\mathbf{f}_{2} +  v^{n-1}_{x}u^{n-1}+v^{n-1}_{y}v^{n-1}.
			\end{split}
		\end{align}
	\end{subequations}
	\item {\bf Iterative method III}
	\begin{subequations}
		\begin{align}
			-\nu \Delta u^{n}  + u^{n}_{x}u^{n-1}+u^{n}_{y}v^{n-1} + p_{x} &=\mathbf{f}_{1},& & \\
			-\nu \Delta v^{n}  +  v^{n}_{x}u^{n-1}+v^{n}_{y}v^{n-1} + p_{y} &=\mathbf{f}_{2}.& &
		\end{align}
	\end{subequations}
	
\end{itemize}

\section{Linearized Learning algorithm with Multiscale Deep Neural Network}
\label{sec:Multiscale}

\subsection{Multiscale deep neural network}
\begin{figure}[ptbh]
\centering
\includegraphics[width=0.65\linewidth]{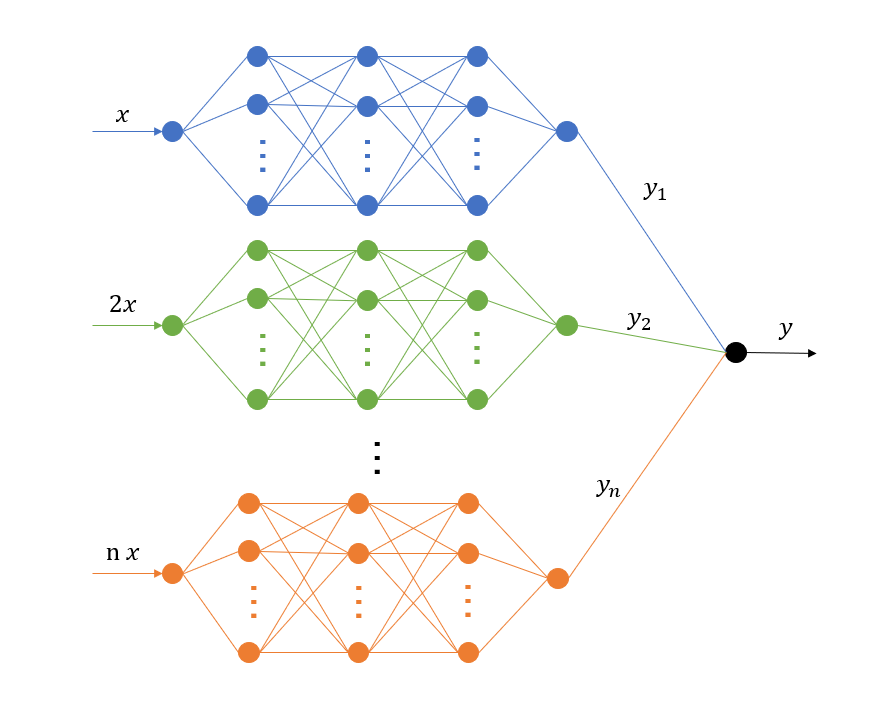}\caption{Illustration of a
MscaleDNN.}
\label{net}
\end{figure}

In order to improve the capability of the DNN to represent functions with multiple scales, we will apply the MscaleDNN \cite{CiCP-28-1970}, which consists of a series of parallel normal sub-neural networks. Each of the sub-networks will receive a scaled version of the input and their outputs will then be combined to make the final output of the MscaleDNN (refer to Fig. \ref{net}). The individual sub-network in the MscaleDNN with a scaled input is designed to approximate a segment of frequency content of the targeted function and the effect of the scaling is to convert a specific high frequency content to a lower frequency range so the learning can be accomplished much quickly. Due to the radial scaling used in the MscaleDNN, it is specially fitting for approximation of high dimensional PDEs.

Fig. \ref{net} shows the schematics of a MscaleDNN consisting of $n$ networks.
Each scaled input passing through a sub-network of $L$ hidden layers can be
expressed in the following formula
\begin{equation}
f_{\vtheta}(\vx)=\vW^{[L-1]}\sigma\circ(\cdots(\mW^{[1]}\sigma\circ
(\mW^{[0]}(\vx)+\vb^{[0]})+\vb^{[1]})\cdots)+\vb^{[L-1]},\label{dnn}%
\end{equation}
where $W^{[1]} \cdots W^{[L-1]}$ $b^{[1]} \cdots b^{[L-1]}$ are the weight
matrices and bias unknowns, respectively, to be optimized via the training,
$\sigma(x)$ is the activation function. In this work, the following plane wave
activation function will be used due to its localized frequency property
\cite{Wang_2020},
\begin{equation}
\sigma(x)=\sin(x).
\end{equation}

For the input scales, we could select the scale for the $i$-th sub-network to
be $i$ (as shown in Fig. \ref{net}) or $2^{i-1}$. Mathematically, a MscaleDNN
solution $f(\vx)$ is represented by the following sum of sub-networks
$f_{\theta^{n_{i}}}$ with network parameters denoted by $\theta^{n_{i}}$
(i.e., weight matrices and bias)
\begin{equation}
f(\vx)\sim{\displaystyle\sum\limits_{i=1}^{M}}f_{\theta^{n_{i}}}(\alpha
_{i}\vx), \label{f_app}%
\end{equation}
where $\alpha_{i}$ is the chosen scale for the $i$-th sub-network in Fig.
\ref{net}. For more details on the design and discussion of the MscaleDNN,
please refer to \cite{CiCP-28-1970}.

\subsection{Four Linearization Schemes}
\label{sec:linear}
To speed up the convergence of the training, we propose an iterative solution procedure for the stationary Navier Stokes equation so that each step we train a linearized version of the Navier-Stokes equation \cite{HE20091351}.

Assume the current learned velocities are \(\left( u^{tmp}, v^{tmp}\right)\) and the new velocities to be learnt are denoted as \(\left( u^{},v^{} \right)\), then we can linearize the Navier-Stokes equation in the following four schemes,
\begin{itemize}
    \item {\bf Scheme 1:} GradFixed: Gradients of velocities in the nonlinear term will be fixed during training and updated after some specific epochs.
 \begin{subequations}
        \begin{align}
        -\nu \Delta u^{}  + u^{tmp}_{x}u^{}+u^{tmp}_{y}v^{} + p_{x} &=\mathbf{f}_{1},& & \\
            -\nu \Delta v^{}  +  v^{tmp}_{x}u^{}+v^{tmp}_{y}v^{} + p_{y} &=\mathbf{f}_{2},& & \\
            \Delta p + 2(-u_{x}^{}v_{y}^{}+ u_{y}^{}v_{x}^{})&= \nabla\cdot \mathbf{f}. & &
        \end{align}
        \label{loss:GradFixed}
\end{subequations}
    \item {\bf Scheme 2:} VFixed: Velocities in the nonlinear term are fixed during training and updated after every specific epochs.
\begin{subequations}
    \begin{align}
        -\nu \Delta u^{}  + u^{}_{x}u^{tmp}+u^{}_{y}v^{tmp} + p_{x} &=\mathbf{f}_{1},& & \\
            -\nu \Delta v^{}  +  v^{}_{x}u^{tmp}+v^{}_{y}v^{tmp} + p_{y} &=\mathbf{f}_{2},& & \\
            \Delta p + 2(-u_{x}^{}v_{y}^{}+ u_{y}^{}v_{x}^{})&= \nabla\cdot \mathbf{f}. & &
        \end{align}
        \label{loss:VFixed}
\end{subequations}
      \item {\bf Scheme 3:} VFixed1: This is the strong form of the Iterative method II, which can be seen as A correction of VFixed.
	\begin{subequations}
		\begin{align}
			\begin{split}
				-\nu \Delta u^{}  + (u^{}_{x}u^{tmp}+u^{}_{y}v^{tmp} &+ u^{tmp}_{x}u^{}+u^{tmp}_{y}v^{} )+ p_{x} \\
				&=\mathbf{f}_{1} + u^{tmp}_{x}u^{tmp}+u^{tmp}_{y}v^{tmp},
			\end{split}
		\end{align}
		\begin{align}
			\begin{split}
				-\nu \Delta v^{}  + ( v^{}_{x}u^{tmp}+v^{}_{y}v^{tmp}&+v^{tmp}_{x}u^{}+v^{tmp}_{y}v^{} )+ p_{y} \\
				&=\mathbf{f}_{2} +  v^{tmp}_{x}u^{tmp}+v^{tmp}_{y}v^{tmp}.
			\end{split}
		\end{align}
        \label{loss:VFixed1}
\end{subequations}
	
\item {\bf Scheme 4:} Hybrid: The nonlinear term is replaced by the averages of the above two schemes.

\begin{subequations}
    \begin{align}
        -\nu \Delta u^{}  + \frac{1}{2}(u^{}_{x}u^{tmp}+u^{}_{y}v^{tmp} + u^{tmp}_{x}u^{}+u^{tmp}_{y}v^{} )+ p_{x} &=\mathbf{f}_{1},& & \\
        -\nu \Delta v^{}  + \frac{1}{2}( v^{}_{x}u^{tmp}+v^{}_{y}v^{tmp}+v^{tmp}_{x}u^{}+v^{tmp}_{y}v^{} )+ p_{y} &=\mathbf{f}_{2},& & \\
            \Delta p + 2(-u_{x}^{}v_{y}^{}+ u_{y}^{}v_{x}^{})&= \nabla\cdot \mathbf{f}. & &
        \end{align}
        \label{loss:Mixed}
\end{subequations}
\end{itemize}

It should be noted that the number of epochs of training a given linearized scheme before updating the fixed term will be a hyperparameter to be adjusted carefully.

Based on equations (\ref{loss:GradFixed}),(\ref{loss:VFixed}),(\ref{loss:Mixed}), we can design four loss functions as follows.

\noindent{\bf Loss function for Scheme 1:}
\begin{equation}
\begin{aligned}
    \text{L}&_{\nabla}  = \text{R}_{u}+\text{R}_{v} + \alpha \text{B}_{\boldsymbol{u}} +\beta\text{R}_{p}, \\
    &\text{R}_{u} = \int_{\Omega} \left(  -\nu \Delta {u}^{}  + {u}^{tmp}_{x}{u}^{}+u^{tmp}_{y}v^{} + p_{x} -\mathbf{f}_{1}\right)^2dxdy,\\
    &\text{R}_{v} = \int_{\Omega} \left(      -\nu \Delta v^{}  +  v^{tmp}_{x}u^{}+v^{tmp}_{y}v^{} + p_{y} -\mathbf{f}_{2} \right)^2dxdy,\\
    &\text{R}_{p}  = \int_{\Omega}\left( \Delta p + 2(-u_{x}^{}v_{y}^{}+ u_{y}^{}v_{x}^{}) - \nabla\cdot \mathbf{f} \right)^2dxdy,\\
    &\text{B}_{\boldsymbol{u}} = \int_{\partial \Omega}\left( \boldsymbol{u}^{} - \boldsymbol{u}_{0} \right)^2dS.
\end{aligned}
\label{loss_1}
\end{equation}

\bigskip
\noindent{\bf Loss function for Scheme 2:}
\begin{equation}
\begin{aligned}
    \text{L}&_{\boldsymbol{u}}  = \text{R}_{u}+\text{R}_{v} + \alpha\text{B}_{\boldsymbol{u}} +\beta\text{R}_{p}, \\
    &\text{R}_{u} = \int_{\Omega} \left(  -\nu \Delta u^{}  + u^{}_{x}u^{tmp}+u^{}_{y}v^{tmp} + p_{x} -\mathbf{f}_{1}\right)^2dxdy,\\
    &\text{R}_{v} = \int_{\Omega} \left(      -\nu \Delta v^{}  +  v^{}_{x}u^{tmp}+v^{}_{y}v^{tmp} + p_{y} -\mathbf{f}_{2} \right)^2dxdy,\\
    &\text{R}_{p}  = \int_{\Omega}\left( \Delta p + 2(-u_{x}^{}v_{y}^{}+ u_{y}^{}v_{x}^{}) - \nabla\cdot \mathbf{f} \right)^2dxdy,\\
    &\text{B}_{\boldsymbol{u}} = \int_{\partial \Omega}\left( \boldsymbol{u}^{} - \boldsymbol{u}_{0} \right)^2dS.
\end{aligned}
\label{loss_2}
\end{equation}

\bigskip
\noindent{\bf Loss function for Scheme 3:}
\begin{equation}
\begin{aligned}
    \text{L}&_{H}  = \text{R}_{u}+\text{R}_{v} + \alpha\text{B}_{\boldsymbol{u}} +\beta\text{R}_{p}, \\
    &\text{R}_{u} = \int_{\Omega} \left(  -\nu \Delta u^{}  +(u^{}_{x}u^{tmp}+u^{}_{y}v^{tmp} + u^{tmp}_{x}u^{}+u^{tmp}_{y}v^{} -u^{tmp}_{x}u^{tmp}-u^{tmp}_{y}v^{tmp} ) + p_{x} -\mathbf{f}_{1}\right)^2dxdy,\\
    &\text{R}_{v} = \int_{\Omega} \left( -\nu \Delta v^{}  +( v^{}_{x}u^{tmp}+v^{}_{y}v^{tmp}+v^{tmp}_{x}u^{}+v^{tmp}_{y}v^{}-v^{tmp}_{x}u^{tmp}-v^{tmp}_{y}v^{tmp}  ) + p_{y} -\mathbf{f}_{2} \right)^2dxdy,\\
    &\text{R}_{p}  = \int_{\Omega}\left( \Delta p + 2(-u_{x}^{}v_{y}^{}+ u_{y}^{}v_{x}^{}) - \nabla\cdot \mathbf{f} \right)^2dxdy,\\
    &\text{B}_{\boldsymbol{u}} = \int_{\partial \Omega}\left( \boldsymbol{u}^{} - \boldsymbol{u}_{0} \right)^2dS.
\end{aligned}
\label{loss_4}
\end{equation}

\bigskip
\noindent{\bf Loss function for Scheme 4:}
\begin{equation}
\begin{aligned}
    \text{L}&_{H}  = \text{R}_{u}+\text{R}_{v} + \alpha\text{B}_{\boldsymbol{u}} +\beta\text{R}_{p}, \\
    &\text{R}_{u} = \int_{\Omega} \left(  -\nu \Delta u^{}  +\frac{1}{2}(u^{}_{x}u^{tmp}+u^{}_{y}v^{tmp} + u^{tmp}_{x}u^{}+u^{tmp}_{y}v^{} ) + p_{x} -\mathbf{f}_{1}\right)^2dxdy,\\
    &\text{R}_{v} = \int_{\Omega} \left( -\nu \Delta v^{}  +\frac{1}{2}( v^{}_{x}u^{tmp}+v^{}_{y}v^{tmp}+v^{tmp}_{x}u^{}+v^{tmp}_{y}v^{} ) + p_{y} -\mathbf{f}_{2} \right)^2dxdy,\\
    &\text{R}_{p}  = \int_{\Omega}\left( \Delta p + 2(-u_{x}^{}v_{y}^{}+ u_{y}^{}v_{x}^{}) - \nabla\cdot \mathbf{f} \right)^2dxdy,\\
    &\text{B}_{\boldsymbol{u}} = \int_{\partial \Omega}\left( \boldsymbol{u}^{} - \boldsymbol{u}_{0} \right)^2dS.
\end{aligned}
\label{loss_3}
\end{equation}

\subsection{Linearized Learning Algorithms}

Our linearized learning algorithms for the Navier-Stokes equation are implemented through the following steps illustrated in {\bf Algorithm 1}.
\IncMargin{1em}
\begin{algorithm}
\SetAlgoLined
\SetKwData{ga}{\(\gamma\)}\SetKwData{ta}{\(\tau\)} \SetKwData{Adam}{Adam Optimizer}
\SetKwFunction{Train}{Train}\SetKwFunction{Loss}{Loss} \SetKwFunction{Update}{Update}
\SetKwInOut{Input}{input}\SetKwInOut{Output}{output}
\Input{\(\boldsymbol{u}^{0},p^{0}, \boldsymbol{u}^{tmp}\)}
\Output{\(\boldsymbol{u}_{},p_{}\)}
\BlankLine
\ta \(=10^{12}\)
\tcp{The thereshold to update the networks}
\ga \(=0.9\)
\tcp{The ratio to make sure the loss is strictly less then the thereshold when updating the networks}
\(\boldsymbol{u}^{tmp}\leftarrow \boldsymbol{u}^{0}\) \;
\(\boldsymbol{u}\leftarrow \boldsymbol{u}^{0}\)\;
\For{$i\leftarrow 1$ \KwTo $N$}{
		\For(\tcp*[h]{\(c\) is a variable to determine the epochs to train the new network  \(\boldsymbol{u}_{}\)}){$j\leftarrow 1$ \KwTo $c$}{\label{forins}
				$L\leftarrow$ \Loss{$\boldsymbol{u}^{tmp},\boldsymbol{u}_{},p_{}$,} \tcp{The \Loss is one of (\ref{loss_1}),(\ref{loss_2}) or (\ref{loss_3})}
			 \Update \(\boldsymbol{u}\) by \Adam with \(L\)\;
			 \Update  \(p\) by \Adam with \(L\)\;
}
\If{\(L\leq\gamma\)\ta}{\label{ut}
\ta \(= L\)\;
\(\boldsymbol{u}_{tmp}\leftarrow \boldsymbol{u}_{ }\)\;
}
}
\caption{disjoint decomposition}\label{algo_disjdecomp}
\end{algorithm}\DecMargin{1em}

\section{Numerical Results}
\label{sec:results}
\subsection{A benchmark problem}

\label{results}
\begin{figure}[ptbh]
\centering
\includegraphics[width=0.75\linewidth]{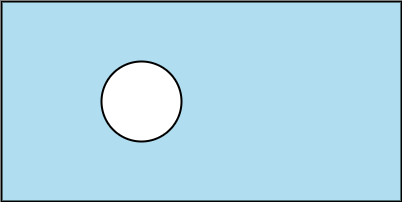}\caption{The Computational Domain.
}
\label{domain}
\end{figure}

We first consider the problem in a rectangle domain \(\Omega = [0,2]\times [0,1]\) with one cylinder hole centered at \(\left( 0.7,0.5 \right)\) whose radius is 0.2 as shown in Figure \ref{domain} , the analytical solutions of the incompressible Navier-Stokes equations is given as follows:
\begin{equation}
    \begin{aligned}
        u&=1-e^{\lambda x} \cos \left(2 m \pi x+2 n \pi y\right), \\
        v&=\frac{\lambda}{2 n \pi} e^{\lambda x} \sin \left(2 m \pi x+2 n \pi y\right)+\frac{m}{n} e^{\lambda x} \cos \left(2 m \pi x+2 n \pi y\right), \\
        p&=\frac{1}{2} (1-e^{2 \lambda x}), \quad \lambda=\frac{\text{Re}}{2}-\sqrt{\frac{\text{Re}^{2}}{4}+4 \pi^{2}}, \quad \text{Re}=\frac{1}{\nu}.
    \end{aligned}
    \label{eq:benchmark_Exact}
\end{equation}

As a benchmark problem, we consider the frequency \(m = 1, n =2\) and the viscosity \(\nu = 0.05\).
The boundary conditions are Dirichlet boundary conditions which are obtained by computing the analytical solutions (\ref{eq:benchmark_Exact}) on boundaries, including the four edges of the ractangle and the circle inside.

\begin{figure}[ptbh]
\centering
\includegraphics[width=0.75\linewidth]{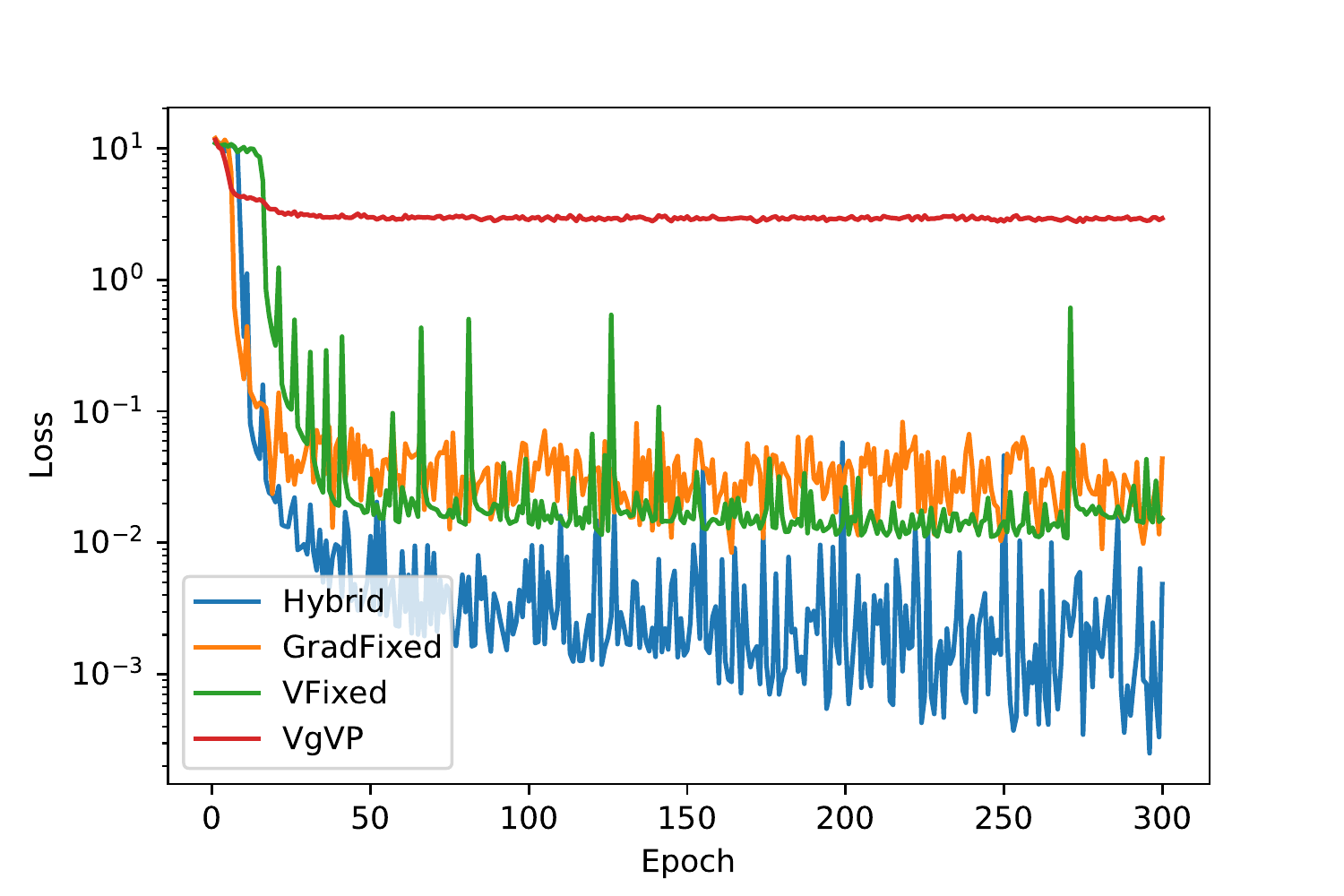}\caption{Losses (bottom 3 lines) of three linearized learning schemes \(\left( \ref{loss_1} \right)-\left( \ref{loss_3} \right)\) and loss (top line) based on nonlinear Navier-Stokes equation (\ref{loss_nl}).
}
\label{loss1to3}
\end{figure}

\begin{figure}
     \centering
     \begin{subfigure}[b]{0.3\textwidth}
         \centering
         \includegraphics[width=\textwidth]{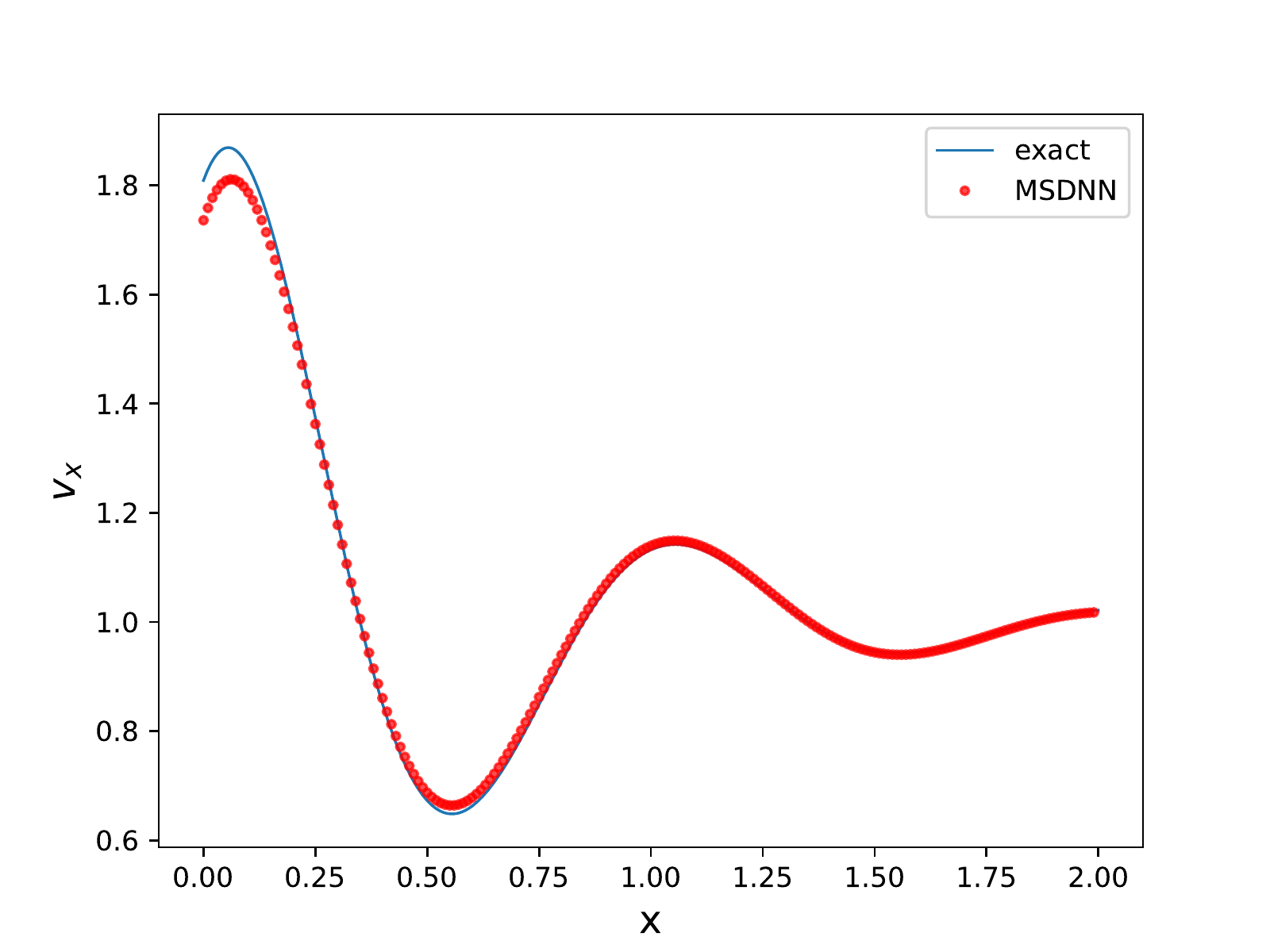}
         \caption{VFixed scheme}
         \label{fig:benchmark_vv}
     \end{subfigure}
     \hfill
     \begin{subfigure}[b]{0.3\textwidth}
         \centering
         \includegraphics[width=\textwidth]{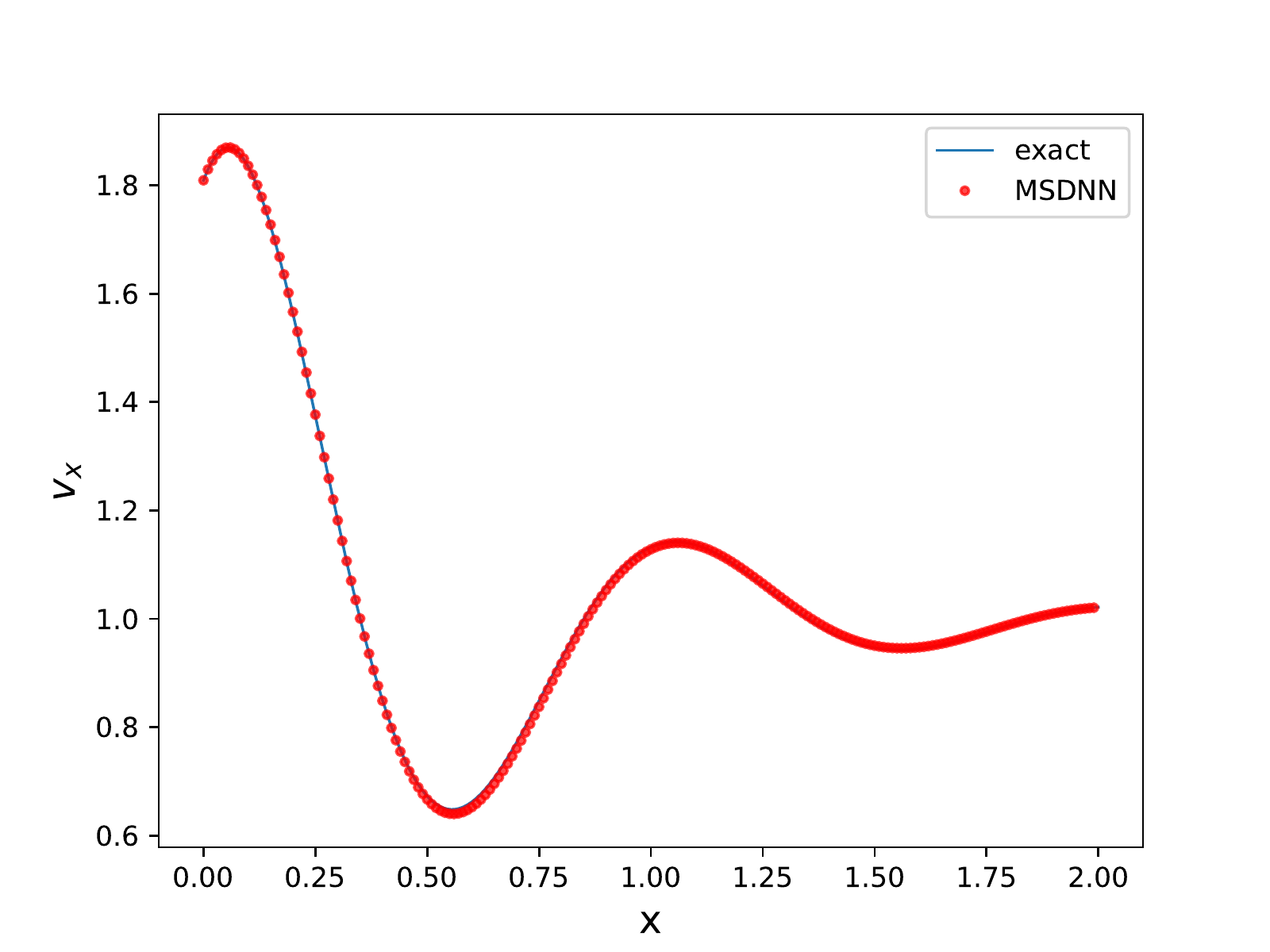}
         \caption{Hybrid scheme}
         \label{fig:benchmark_vmix}
     \end{subfigure}
     \hfill
     \begin{subfigure}[b]{0.3\textwidth}
         \centering
         \includegraphics[width=\textwidth]{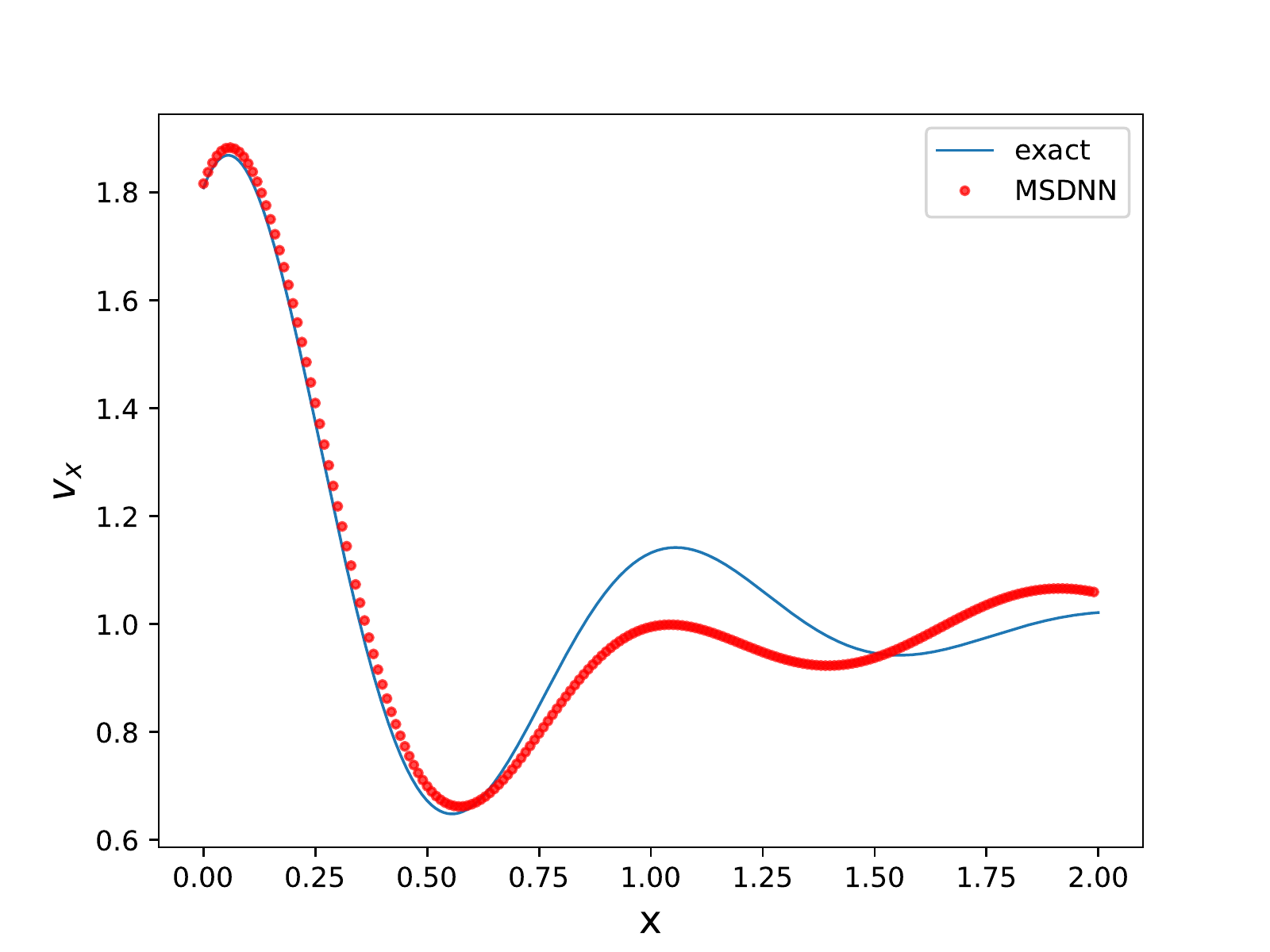}
         \caption{GradFixed scheme}
         \label{fig:benchmark_vgrad}
     \end{subfigure}
        \caption{The x components of velocity after 300 epoch training for bench mark problem along line $y=0.7$}
        \label{fig:benchmarkv}
\end{figure}

\begin{figure}
     \centering
     \begin{subfigure}[b]{0.3\textwidth}
         \centering
         \includegraphics[width=\textwidth]{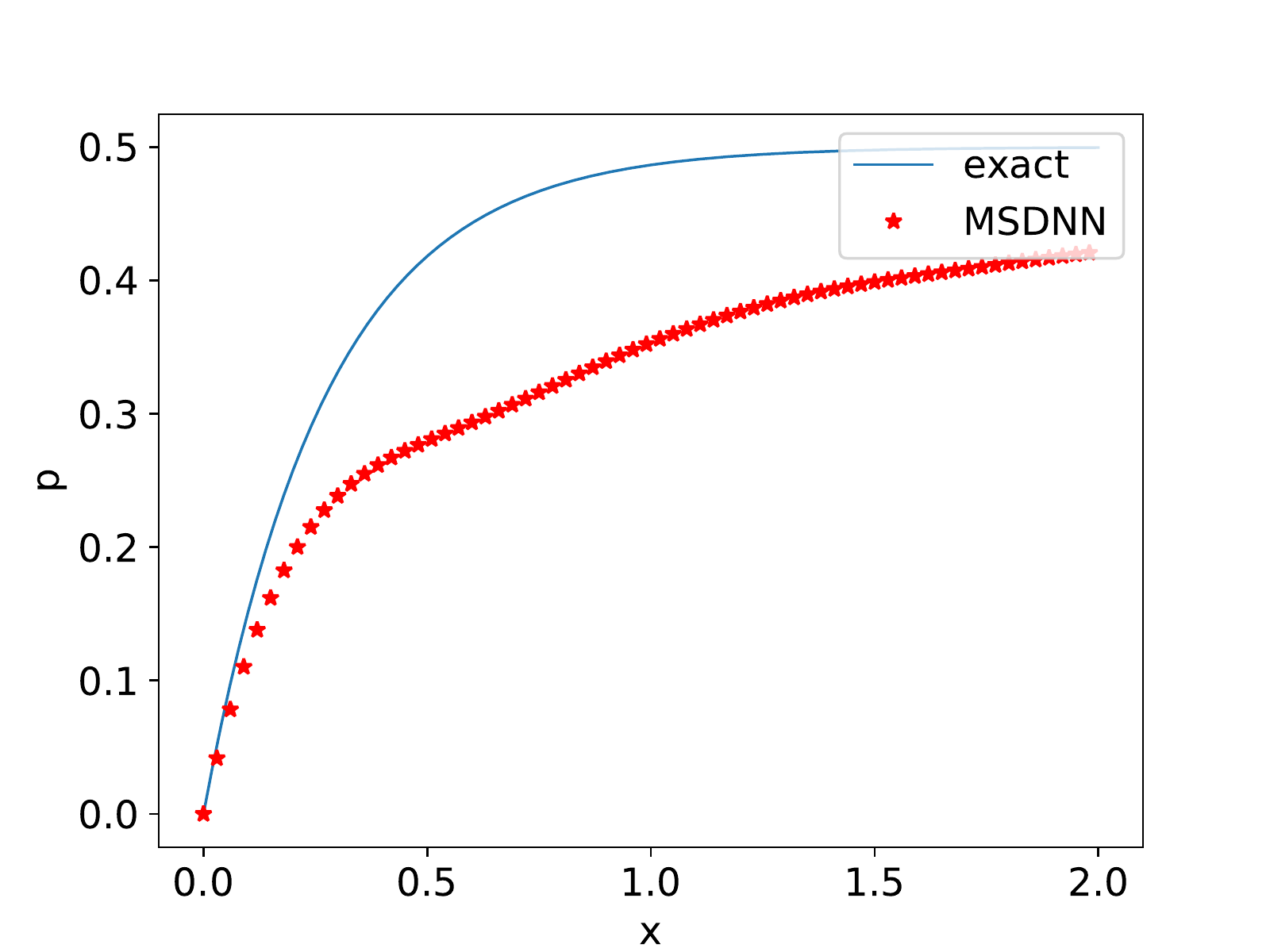}
         \caption{VFixed scheme}
         \label{fig:benchmark_v}
     \end{subfigure}
     \hfill
     \begin{subfigure}[b]{0.3\textwidth}
         \centering
         \includegraphics[width=\textwidth]{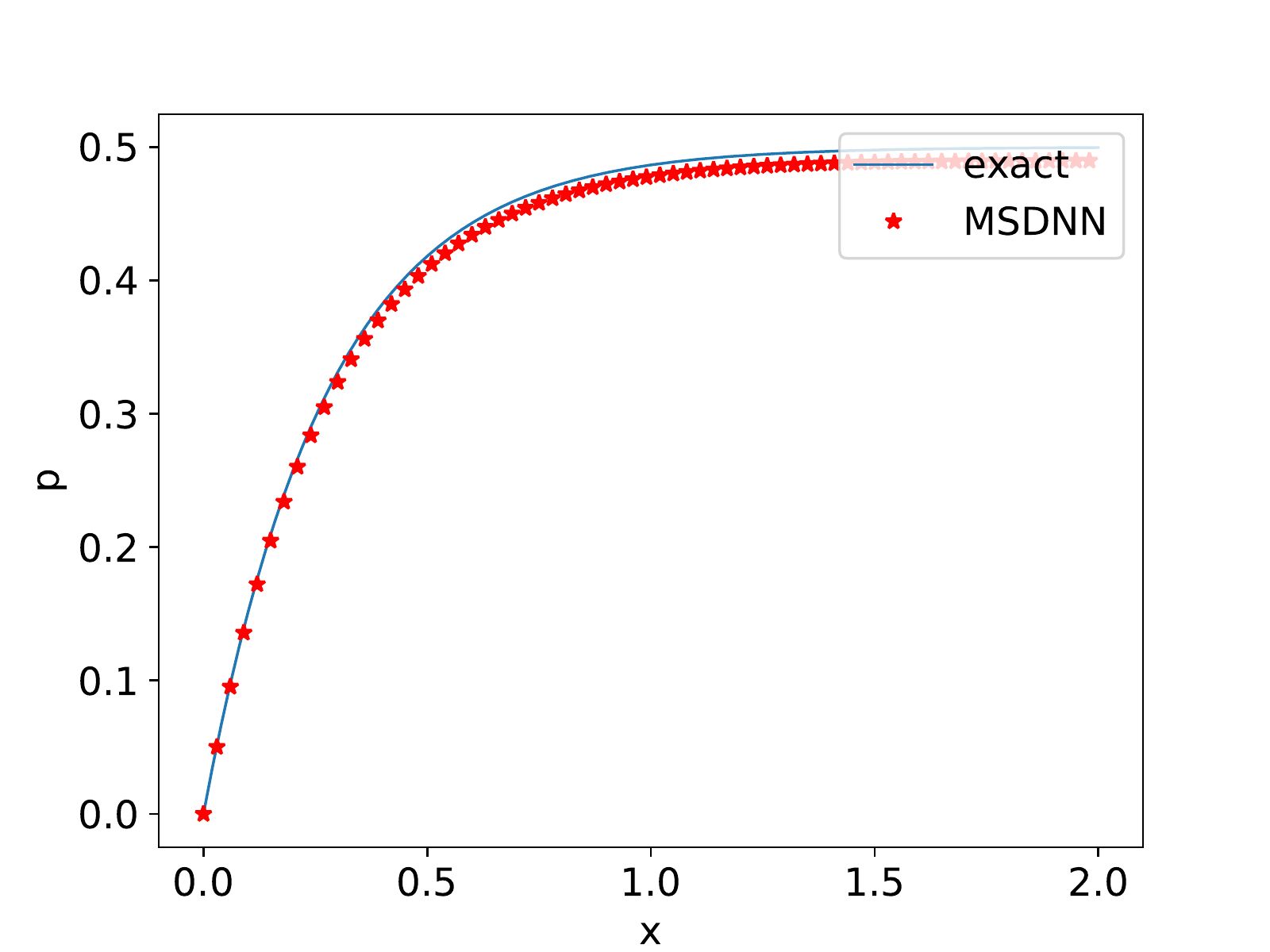}
         \caption{Hybrid scheme}
         \label{fig:benchmark_mix}
     \end{subfigure}
     \hfill
     \begin{subfigure}[b]{0.3\textwidth}
         \centering
         \includegraphics[width=\textwidth]{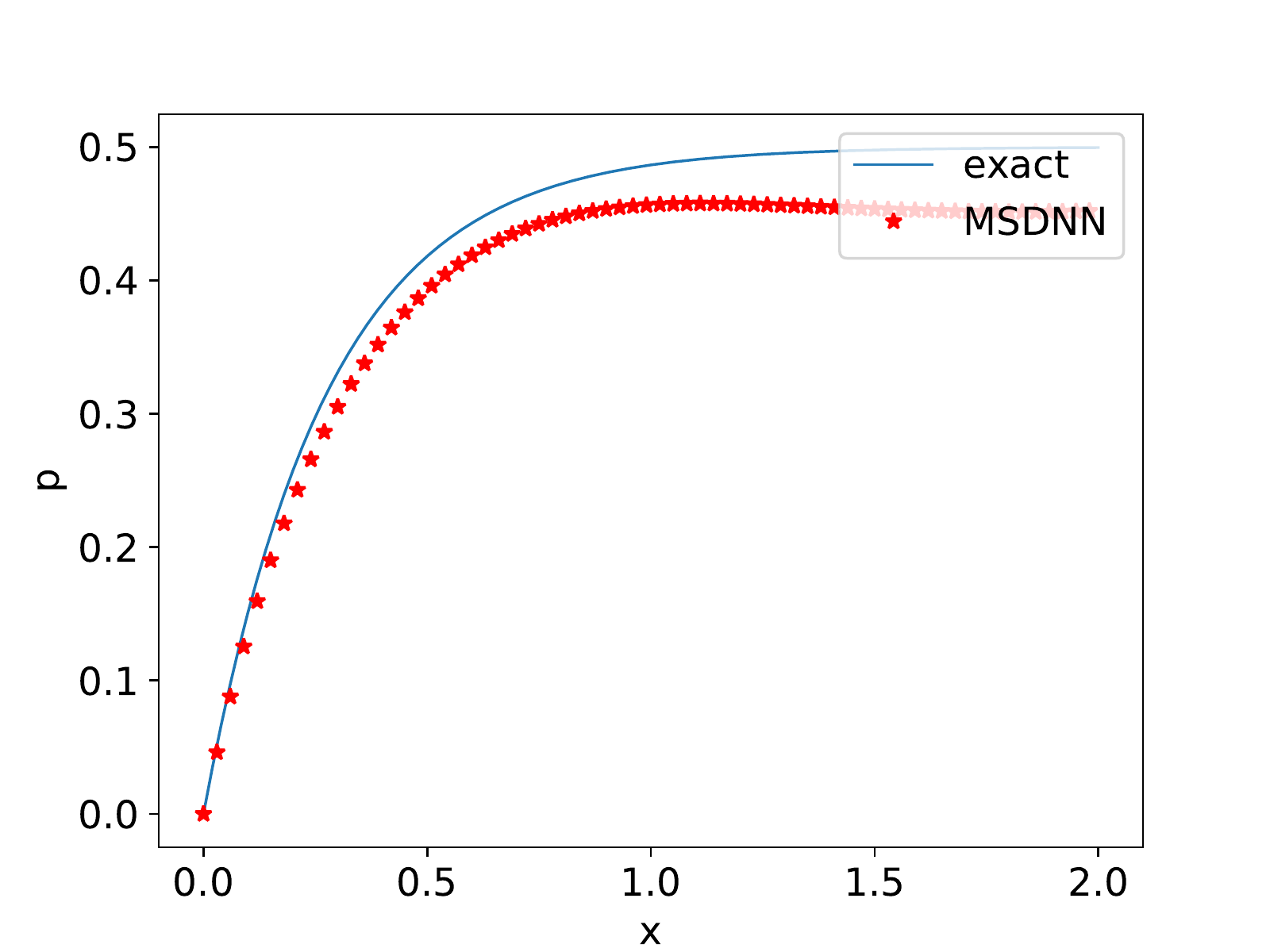}
         \caption{GradFixed scheme}
         \label{fig:benchmark_grad}
     \end{subfigure}
        \caption{The pressures after 300 epoch training for bench mark problem along line $y=0.7$}
        \label{fig:benchmark}
\end{figure}

We compare the performance of the convergence of fcn with different loss schemes, including the VgVP formulation (\ref{loss_nl}) and four linearized schemes (\ref{loss_1}) or (\ref{loss_2}) or (\ref{loss_3}). We randomly sample 160000 points inside \(\Omega\) and 16000 points on \(\partial \Omega\) during each epoch. In the learning process, we set number of batches to be 50 for each epoch. We choose the fully connected neural network with 4 hidden layer, 100 hidden neurons each layer for \(\boldsymbol{u},p\) for all cases. The hyperparameters are the same for four cases. The losses during traing for different cases by minimizing given loss function are compared in Figure \(\ref{loss1to3}\). The results show that all four linearized learning neural networks converge in 300 epochs for all schemes while learning using the loss function (\ref{loss_nl}) for the nonlinear Navier-Stokes equations fails to (top line in Figure \(\ref{loss1to3}\) ).
The comparisons of the x component of velocity and pressure along line \(y = 0.7\) of different linearized schemes after 300 epoch training are shown in Figure \ref{fig:benchmarkv} and \ref{fig:benchmark}.

\subsection{Oscillating flows}
\label{oscillating_simple}
In our previous work \cite{Wang_2020}, it has been shown that a multiscale DNN will be needed to learn oscillatory solutions for the Stokes equations. So the linearized learning schemes with multi-scale deep neural network will be also used for Navier-Stokes equations with oscillatory solutions.  Here we use the same domain like the benchmark problem but the frequencies now are taken to be \(m = 40, n = 35\),  much higher than the benchmark problems and, a training algorithm based on the VGVP scheme using the loss function (\ref{loss_nl}) for the nonlinear NS does not coverage within the same number of epoches. To accelerate the convergence, we adjust the learning rate during training by a decreasing of \(5\%\) every 50 epochs.

\begin{figure}
     \centering
     \begin{subfigure}[b]{0.48\textwidth}
         \centering
         \includegraphics[width=\textwidth]{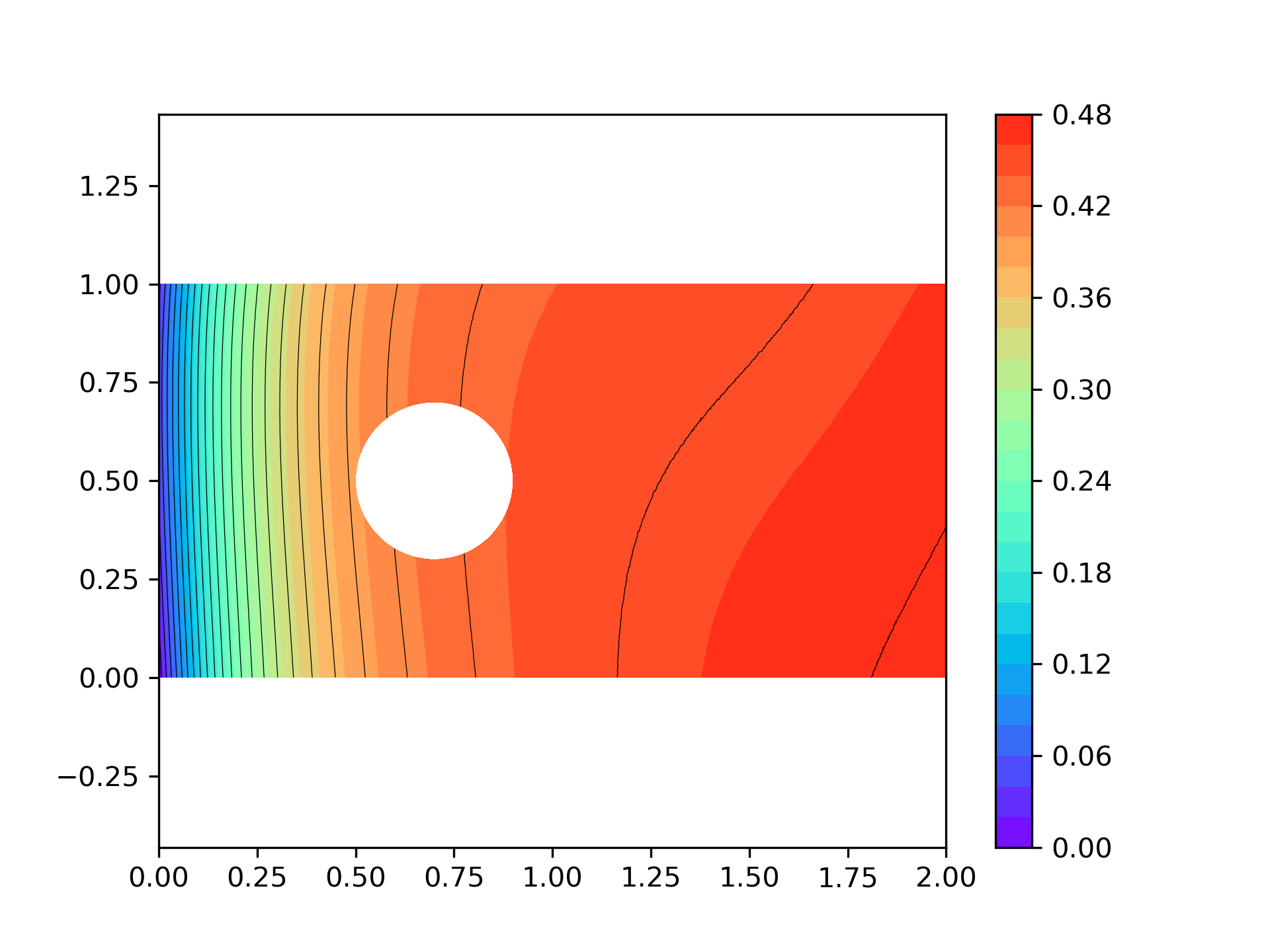}
         \caption{Contour of pressure of the oscillatory case after 1000 epoch training for FCN}
         \label{fig:pressure_Contour_fcn}
     \end{subfigure}
     \hfill
     \begin{subfigure}[b]{0.48\textwidth}
         \centering
         \includegraphics[width=\textwidth]{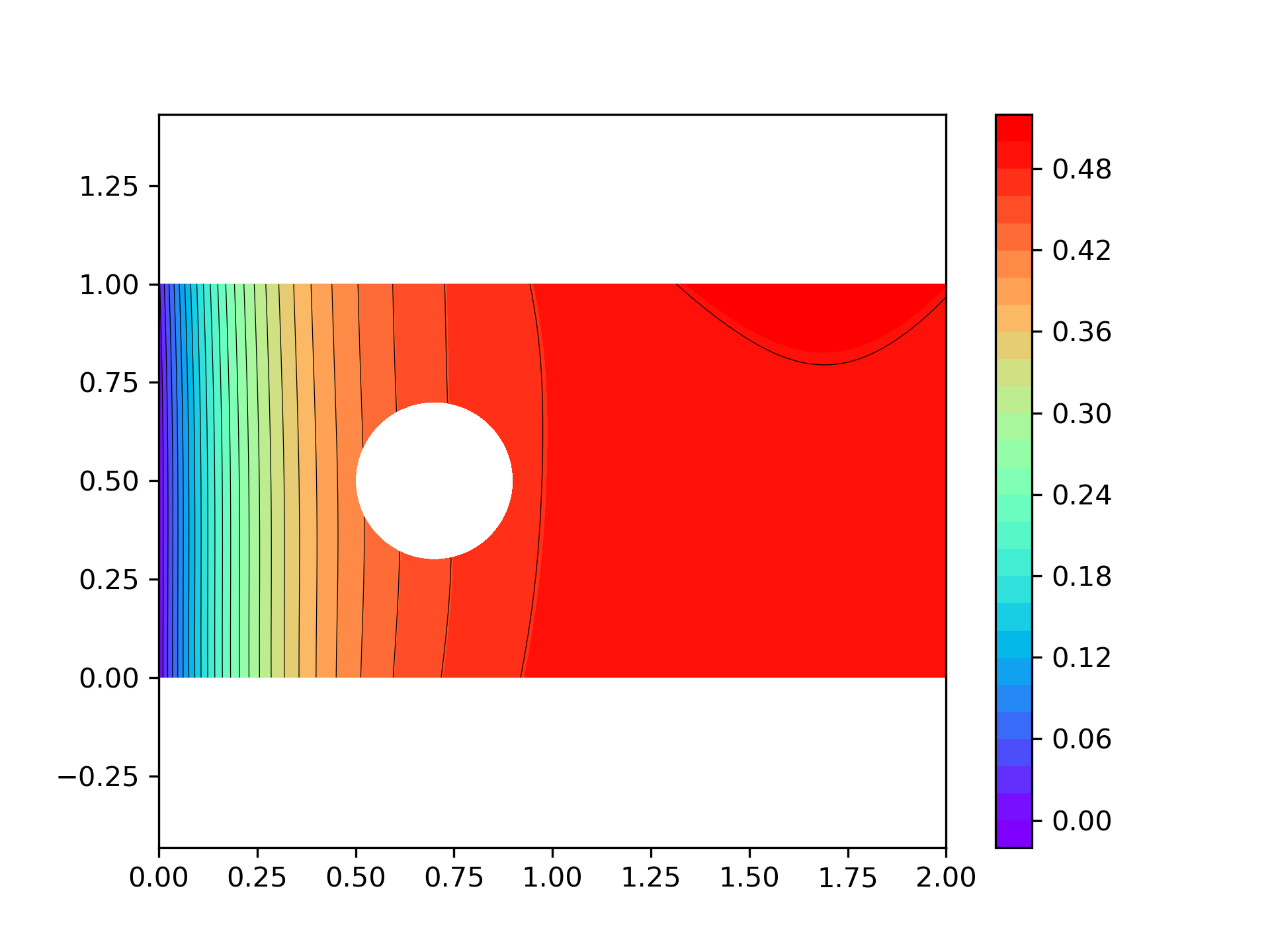}
         \caption{Contour of velocity of the first component after 1000 epoch training for MscaleDNN}
         \label{fig:pressure_Contour_msnn}
     \end{subfigure}
     \hfill
     \caption{The results of pressure of the oscillatory case for different neural networks}
        \label{fig:pressure_Contour}
\end{figure}

\begin{figure}
     \centering
     \begin{subfigure}[b]{0.48\textwidth}
         \centering
         \includegraphics[width=\textwidth]{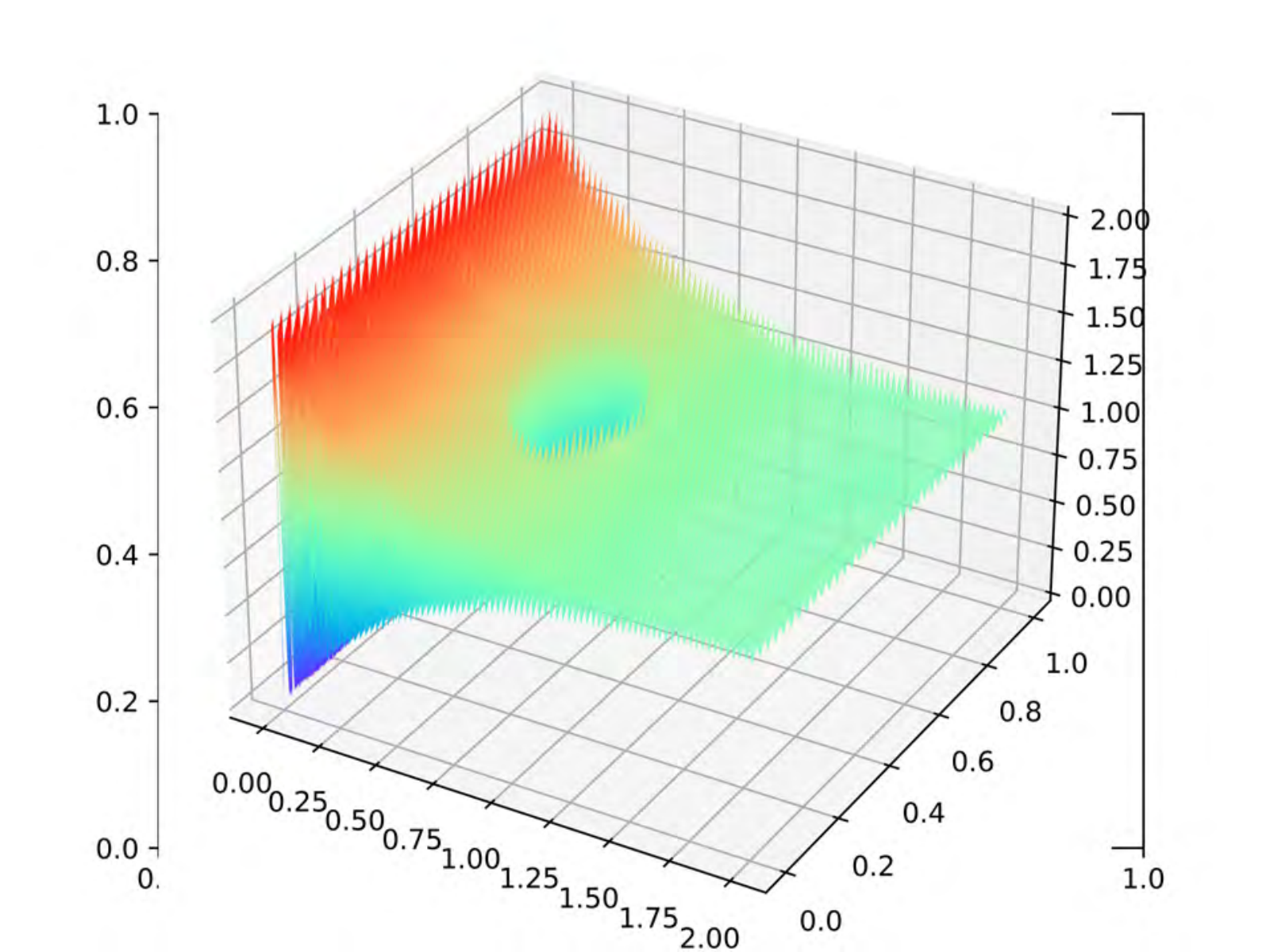}
         \caption{Contour of the first component of velocity of the oscillatory case after 1000 epoch training for FCN}
         \label{fig:velocity_x_Contour_fcn}
     \end{subfigure}
     \hfill
     \begin{subfigure}[b]{0.48\textwidth}
         \centering
         \includegraphics[width=\textwidth]{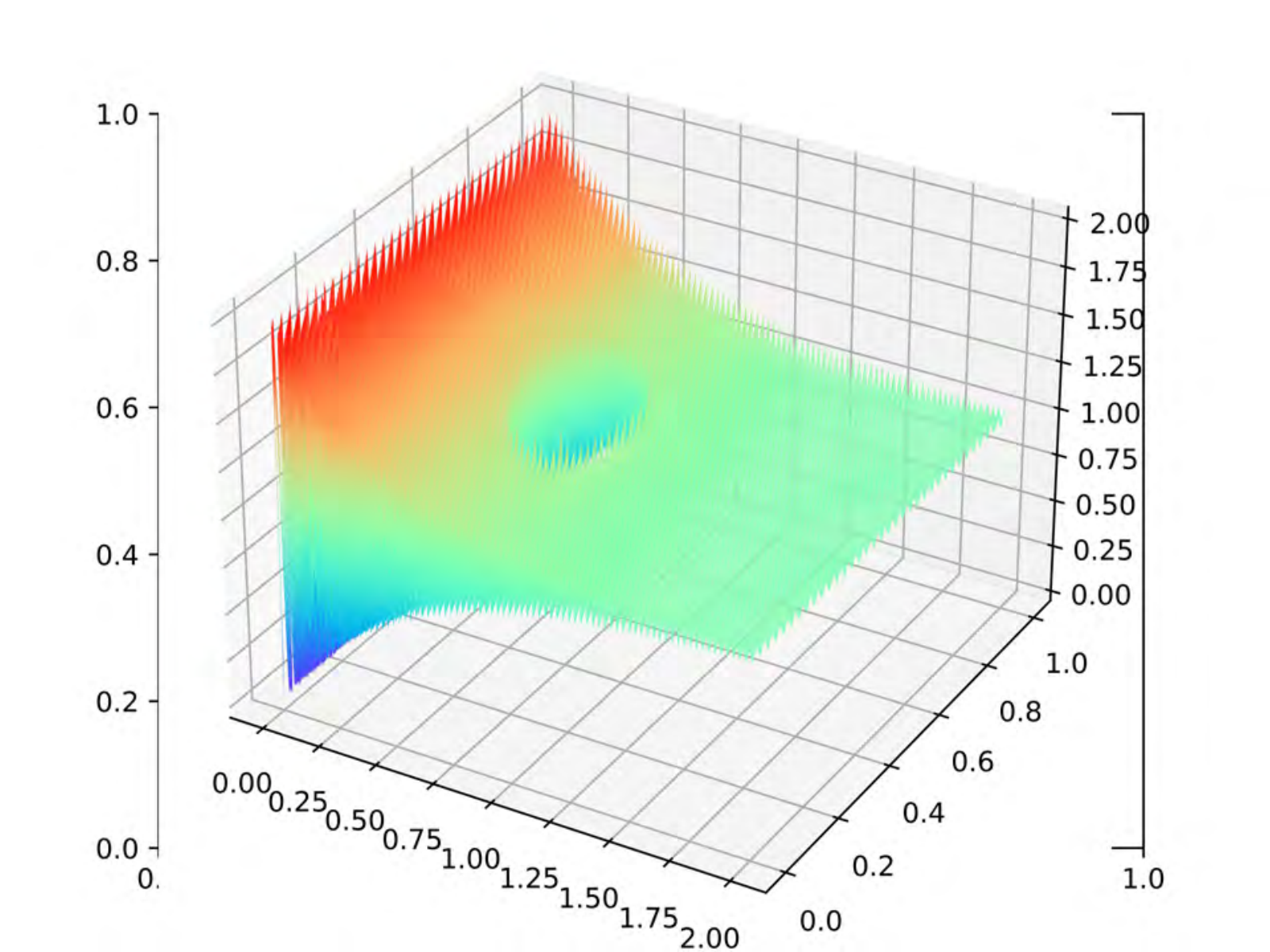}
         \caption{Contour of the first component of velocity of the oscillatory case after 1000 epoch training for MscaleDNN}
         \label{fig:velocity_x_Contour_msnn}
     \end{subfigure}
     \hfill
     \caption{The results of the first component of velocity of the oscillatory case for different neural networks}
        \label{fig:velocity_x_Contour}
\end{figure}

\begin{figure}
     \centering
     \begin{subfigure}[b]{0.48\textwidth}
         \centering
         \includegraphics[width=\textwidth]{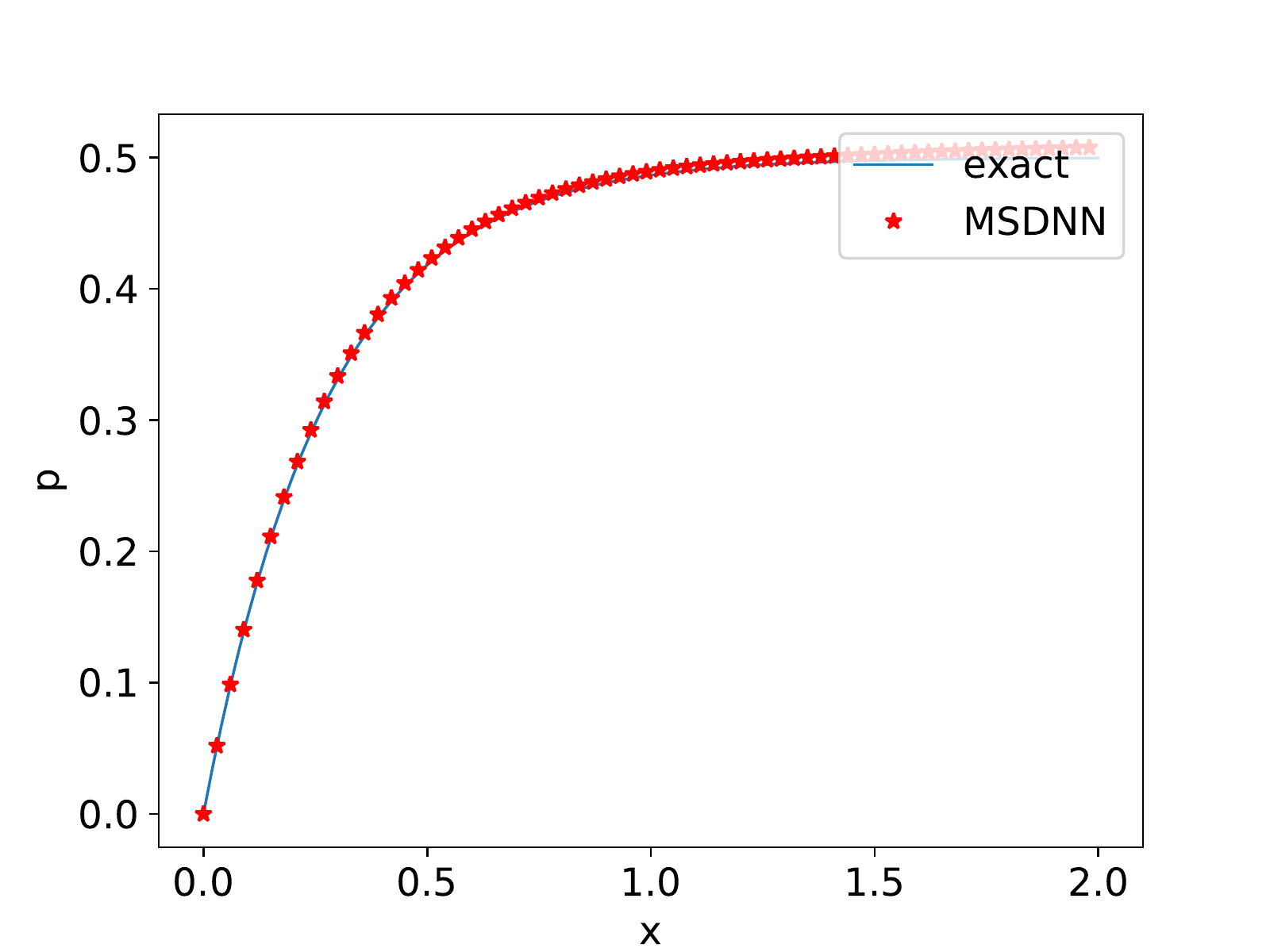}
         \caption{Pressure of the oscillatory case after 1000 epoch training alone line \(y = 0.7\)}
         \label{fig:oscillatory_p}
     \end{subfigure}
     \hfill
     \begin{subfigure}[b]{0.48\textwidth}
         \centering
         \includegraphics[width=\textwidth]{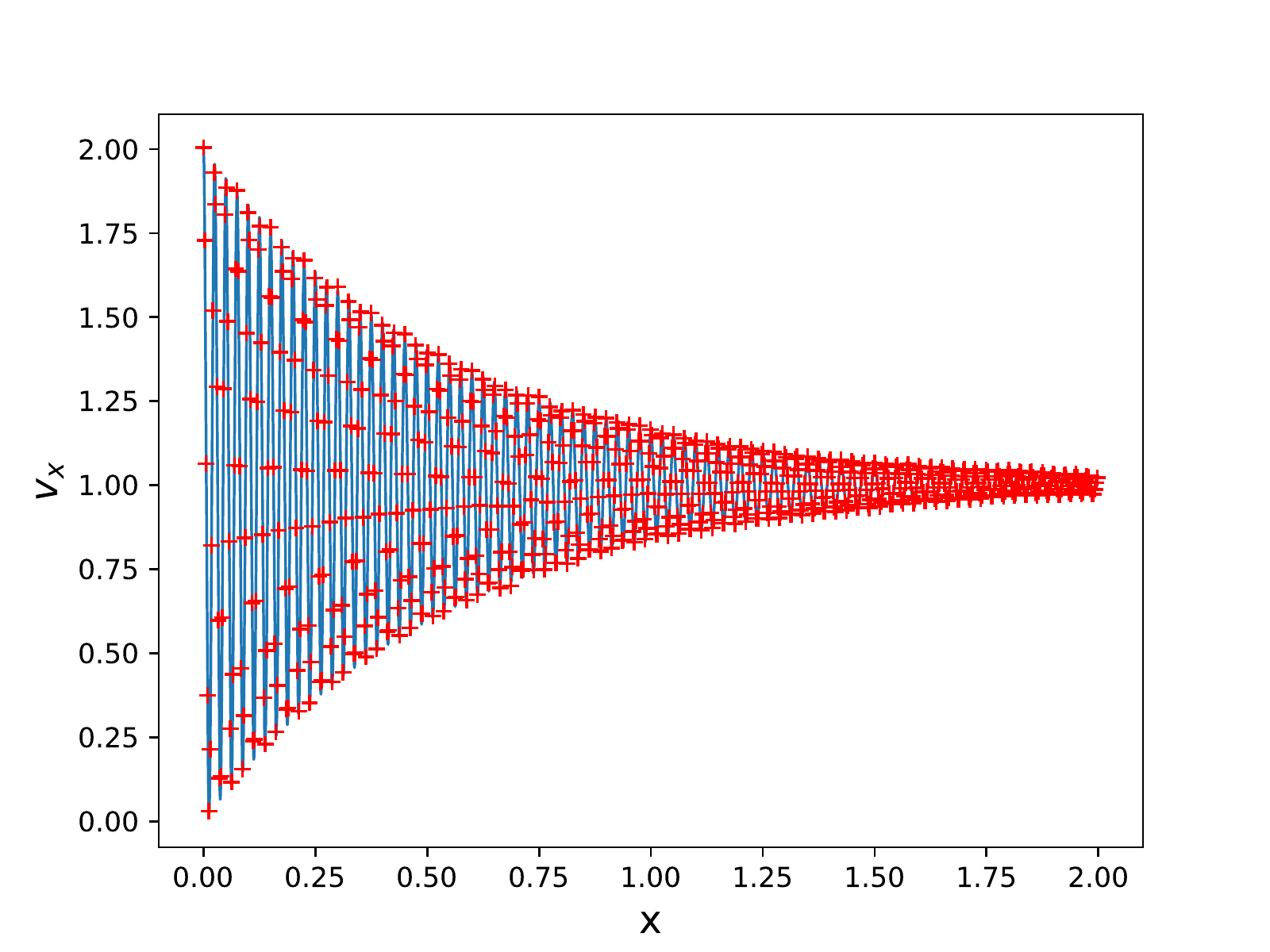}
         \caption{Velocity of the first component after 1000 epoch training along line \(y = 0.7\)}
         \label{fig:oscillatory_u}
     \end{subfigure}
     \hfill
     \caption{The results of the oscillatory case using multi-scale neural networks}
        \label{fig:oscillatory}
\end{figure}

\begin{figure}
     \centering
     \begin{subfigure}[b]{0.48\textwidth}
         \centering
         \includegraphics[width=\textwidth]{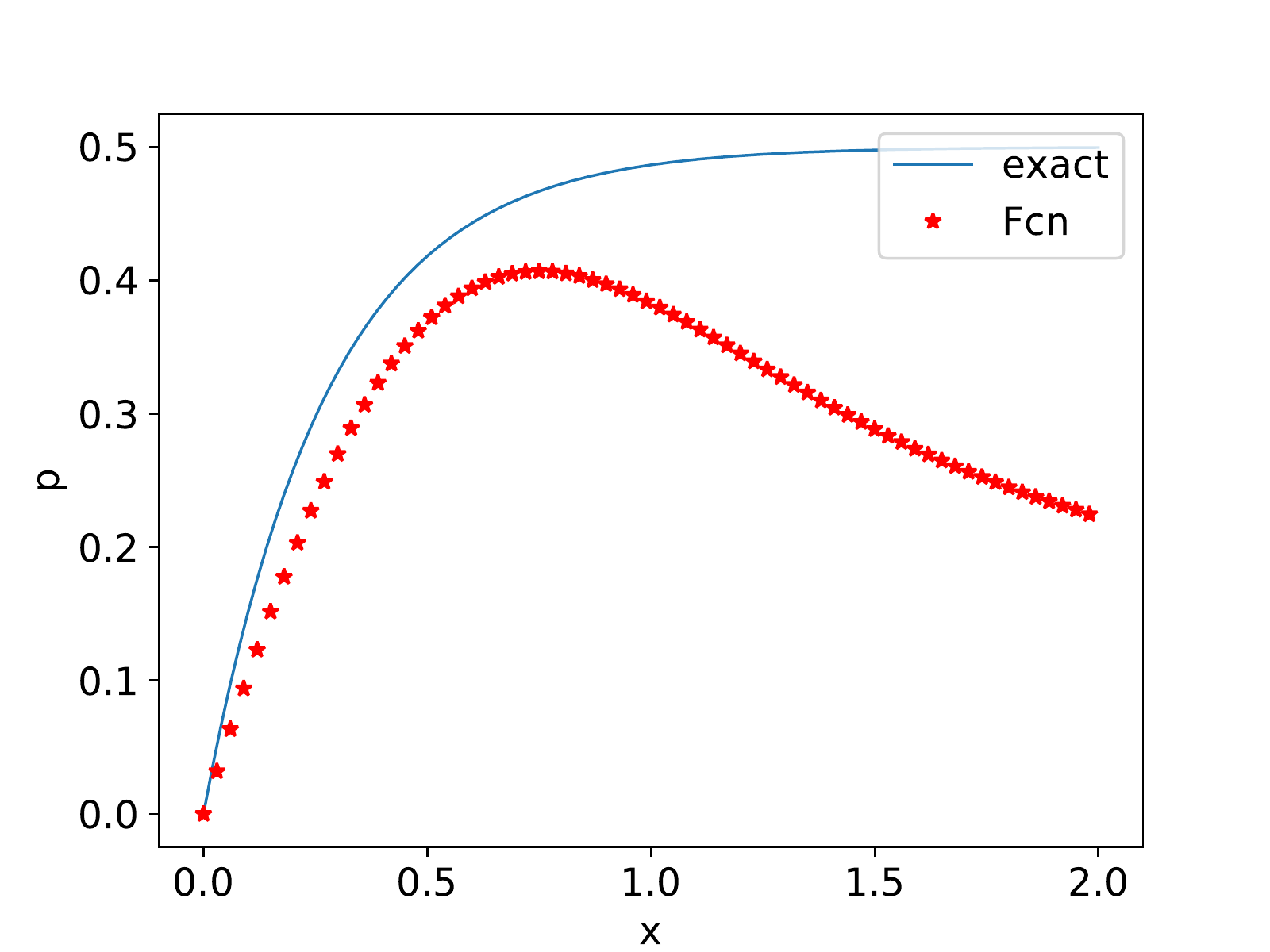}
         \caption{Pressure of the oscillatory case after 1000 epoch training alone line \(y = 0.7\)}
         \label{fig:osci_fcn_p}
     \end{subfigure}
     \hfill
     \begin{subfigure}[b]{0.48\textwidth}
         \centering
         \includegraphics[width=\textwidth]{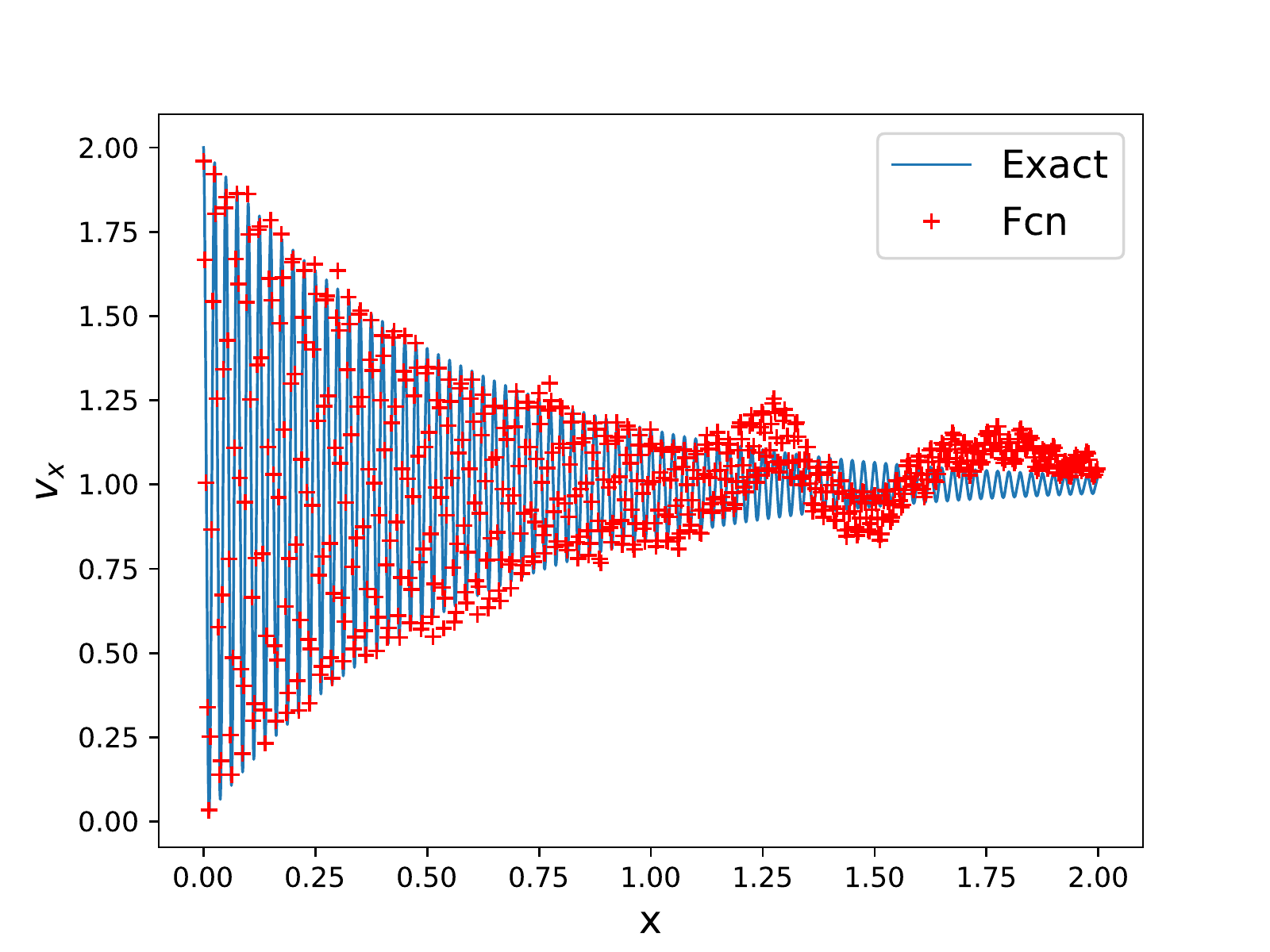}
         \caption{Velocity of the first component after 1000 epoch training along line \(y = 0.7\)}
         \label{fig:osci_fcn_u}
     \end{subfigure}
     \hfill
     \caption{The results of the oscillatory case using fully connected neural networks}
        \label{fig:osci_fcn}
\end{figure}

\begin{figure}
     \centering
     \begin{subfigure}[b]{0.48\textwidth}
         \centering
         \includegraphics[width=\textwidth]{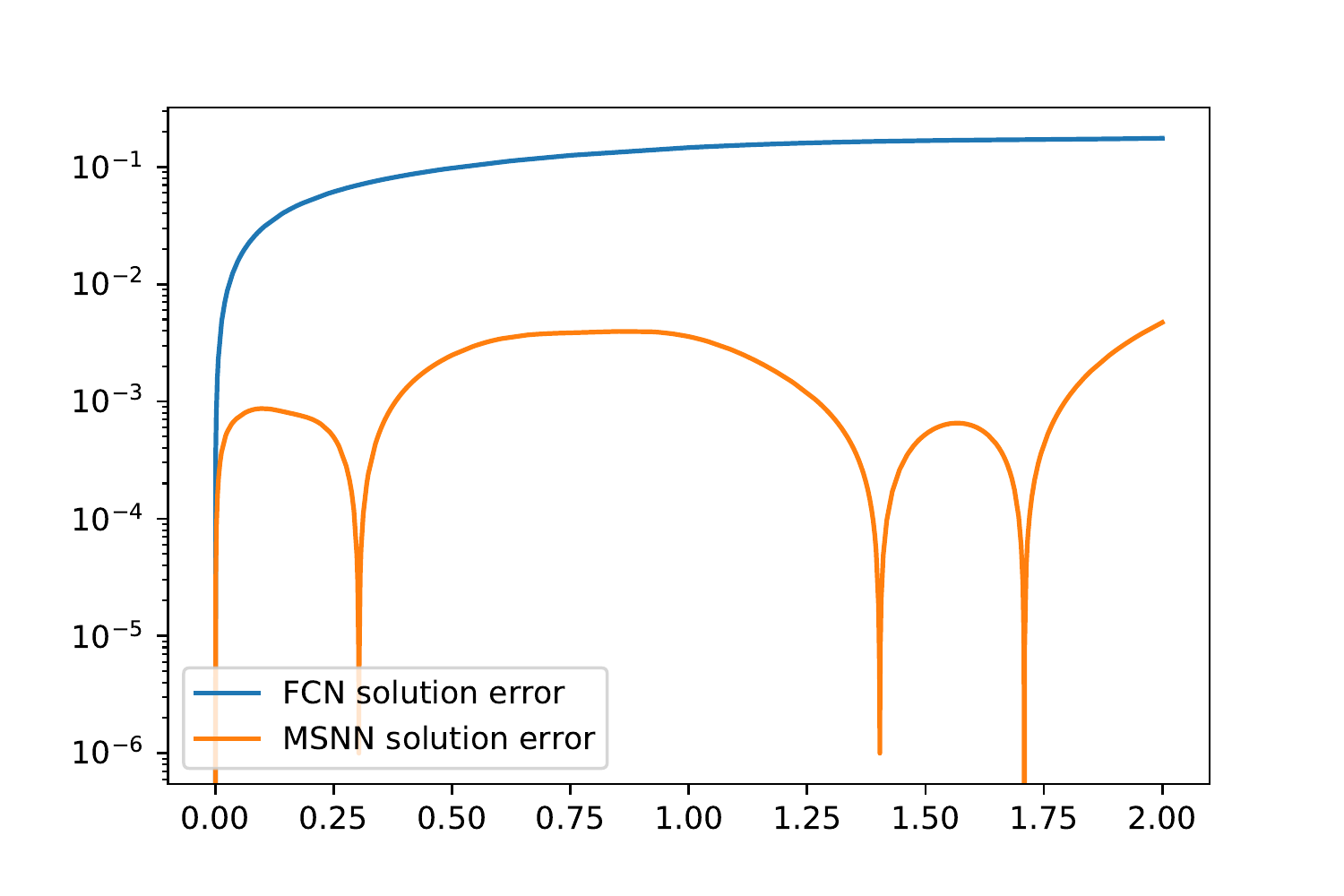}
         \caption{Error of two different models w.r.t. pressure of the oscillatory case after 1000 epoch training alone line \(y = 0.7\)}
         \label{fig:osci_err_p}
     \end{subfigure}
     \hfill
     \begin{subfigure}[b]{0.48\textwidth}
         \centering
         \includegraphics[width=\textwidth]{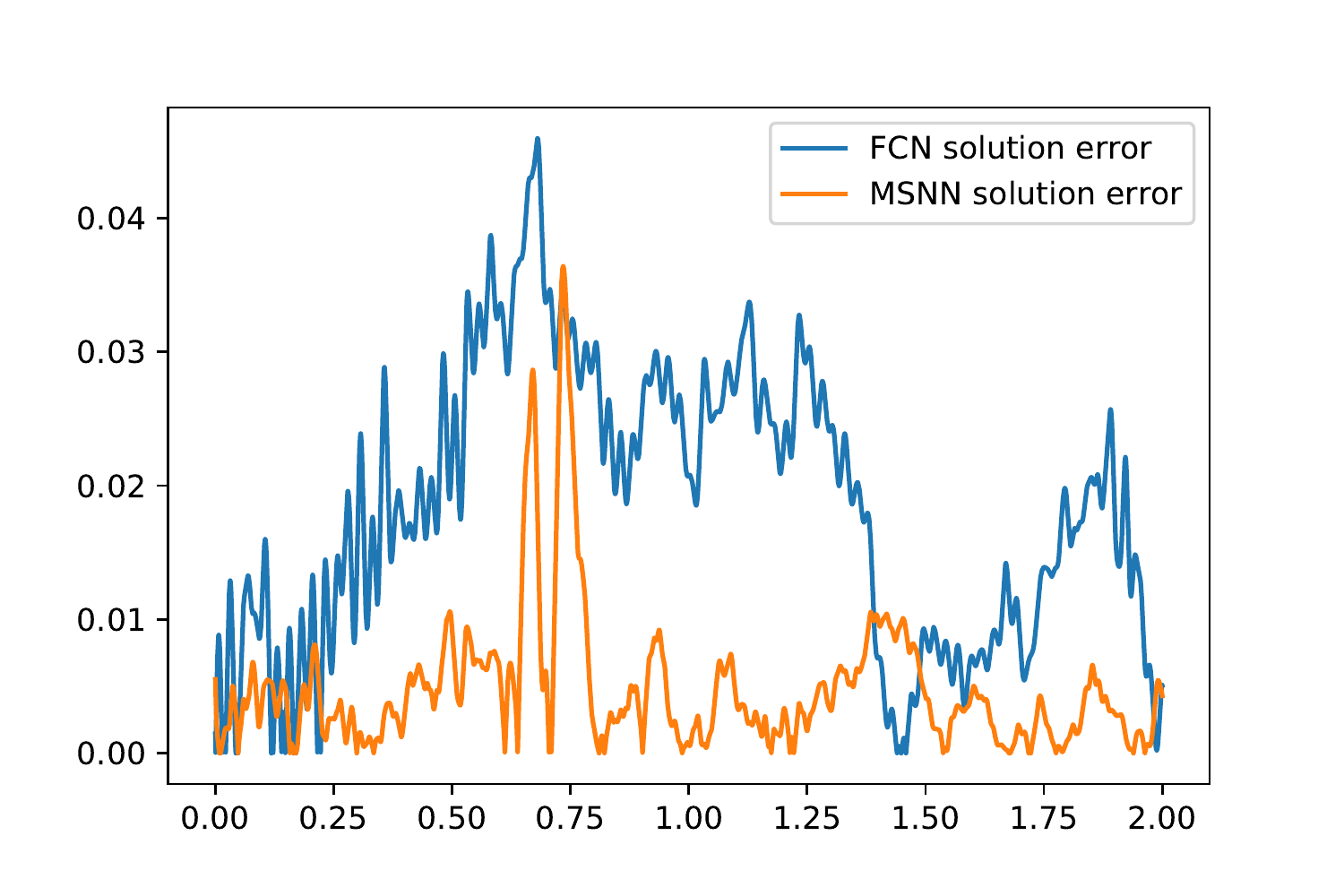}
         \caption{Error of two different models w.r.t. velocity of the first component after 1000 epoch training along line \(y = 0.7\)}
         \label{fig:osci_err_u}
     \end{subfigure}
     \hfill
     \caption{The errors of the oscillatory case for FCN and multiscale network (MSNN) architectures}
        \label{fig:osci_err}
\end{figure}

 We consider Scheme 2 of the linearized learning algorithms where the previous velocity are used to linearize the convection term. The multiscale deep neural networks has 8 scales: \(\{x, 2x, 4x,8x,32x,64x, 128x\}\), whose subnetworks contain 4 hidden layers and 128 hidden neurons in each layer. As a comparison, we also trained a 4-layer fully connected neural network with 1024 hidden neurons combining GradFixed scheme in 1000 epochs. Figures \ref{fig:pressure_Contour} and \ref{fig:velocity_x_Contour} show the predictions of network after 1000 epoch training. Figures \ref{fig:oscillatory} and \ref{fig:osci_fcn} give more details along line \(y = 0.7\).
Fig \ref{fig:osci_err} shows the relative errors of these 2 different neural network structures along line \(y = 0.7\). The MscaleDNN improves the accuracy of both the pressure field and the velocity field compared with the fully connected neural network.



\subsection{Oscillating flows in complex domain}
\begin{figure}[ptbh]
\centering
\includegraphics[width=0.75\linewidth]{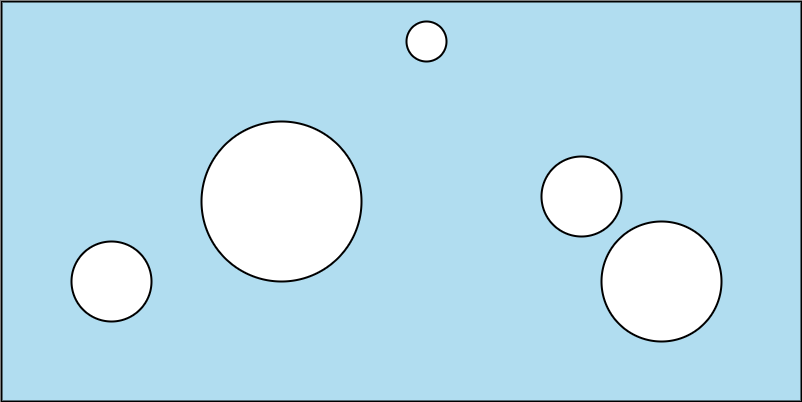}\caption{The more complex domain}
\label{complex_domain}
\end{figure}

 In this section, we consider an oscillating case in complex domain with more than one hole. The domain is shown in Figure \ref{complex_domain}.

We will use scheme vFixed1 \(\left( \ref{loss:VFixed1}\right)\) to compare with the schemes we propose to give an intuitive understanding about the error bounds of scheme GradFixed, VFixed and Hybrid. The rigorous discussion about the error bounds and convergence of these schemes will be studied in the future.

\begin{figure}
     \centering
     \begin{subfigure}[b]{0.48\textwidth}
         \centering
         \includegraphics[width=\textwidth]{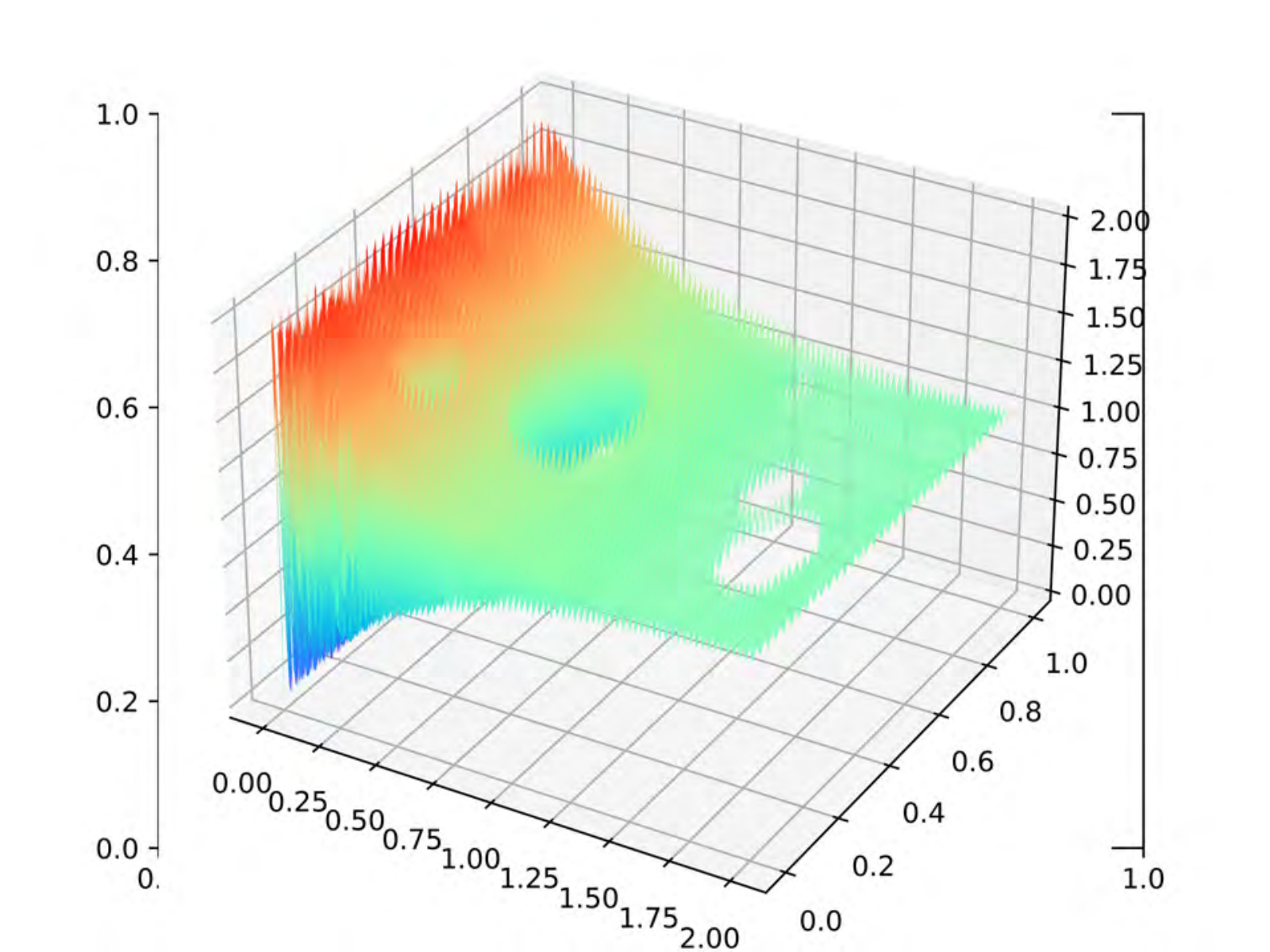}
         \caption{Contour of the first component of velocity of the oscillatory case after 1000 epoch training for scheme vFixed}
         \label{fig:complexDomain_pressure_Contour_vFixed}
     \end{subfigure}
     \hfill
     \begin{subfigure}[b]{0.48\textwidth}
         \centering
         \includegraphics[width=\textwidth]{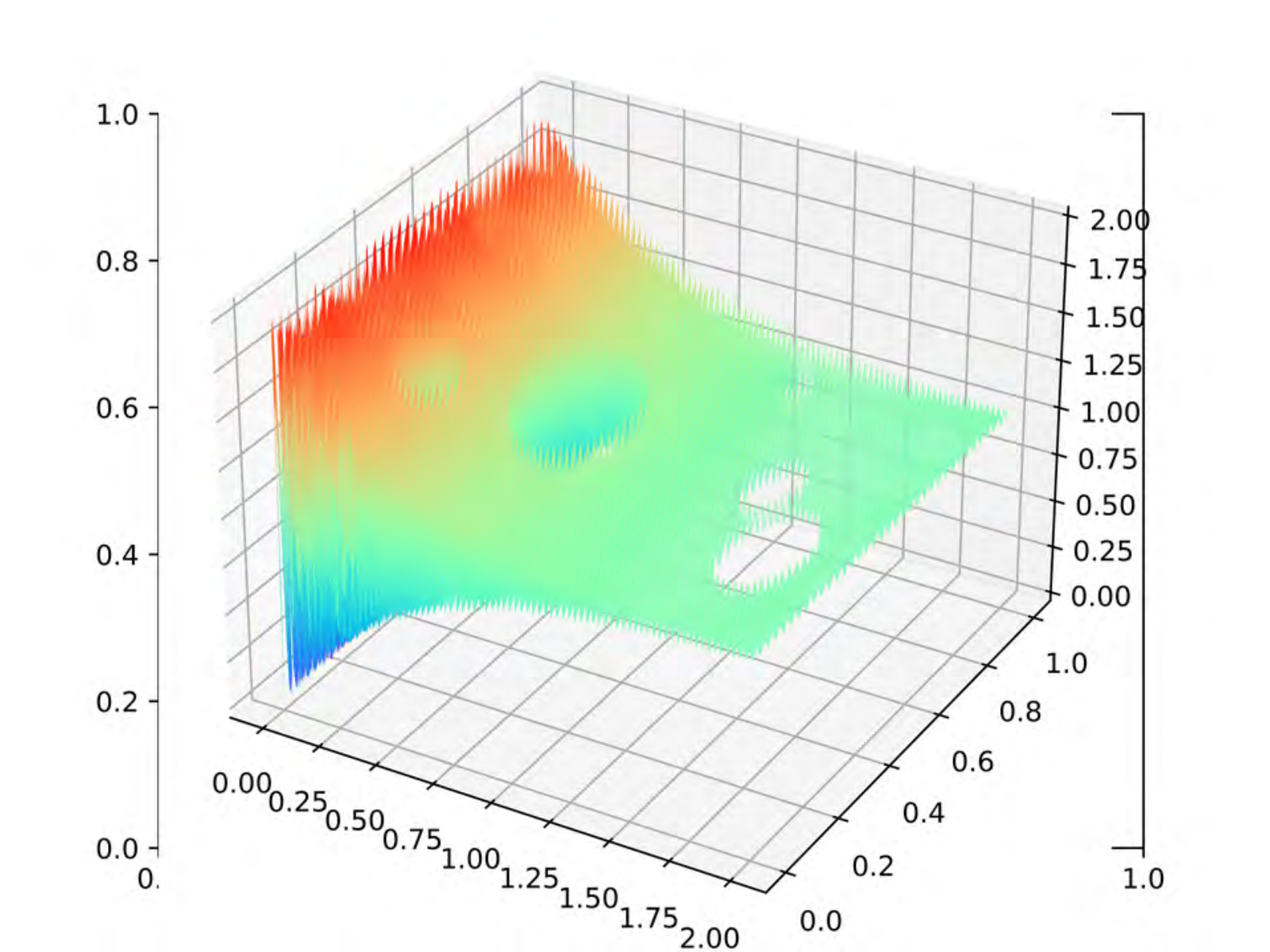}
         \caption{Contour of the first component of velocity of the oscillatory case after 1000 epoch training for scheme vFixed1}
         \label{fig:complexDomain_pressure_Contour_vFixed1}
     \end{subfigure}
     \hfill
     \begin{subfigure}[b]{0.48\textwidth}
         \centering
         \includegraphics[width=\textwidth]{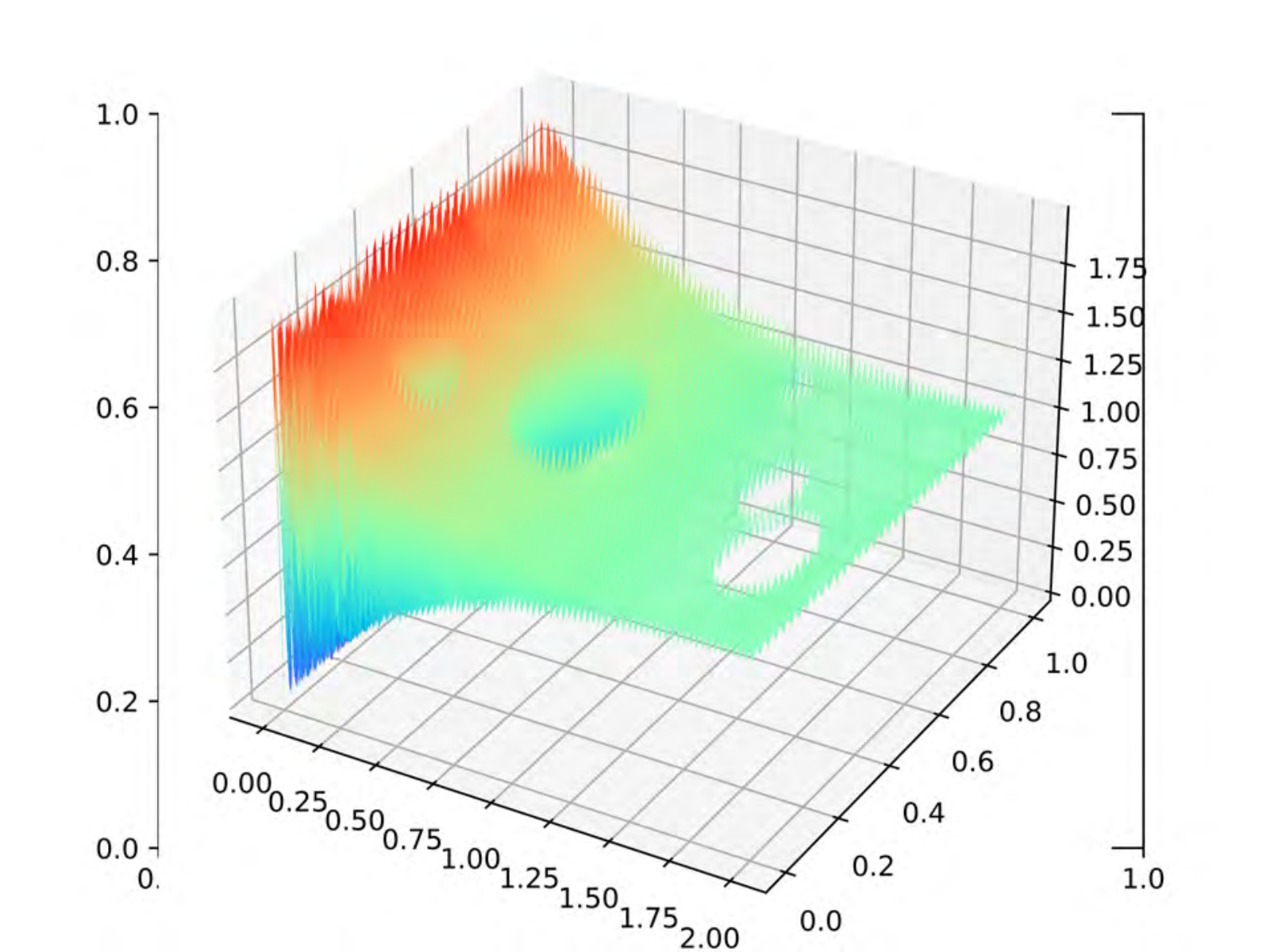}
         \caption{Contour of the first component of velocity of the oscillatory case after 1000 epoch training for scheme gradFixed}
         \label{fig:complexDomain_velocity_x_Contour_gradFixed}
     \end{subfigure}
     \hfill
     \begin{subfigure}[b]{0.48\textwidth}
         \centering
         \includegraphics[width=\textwidth]{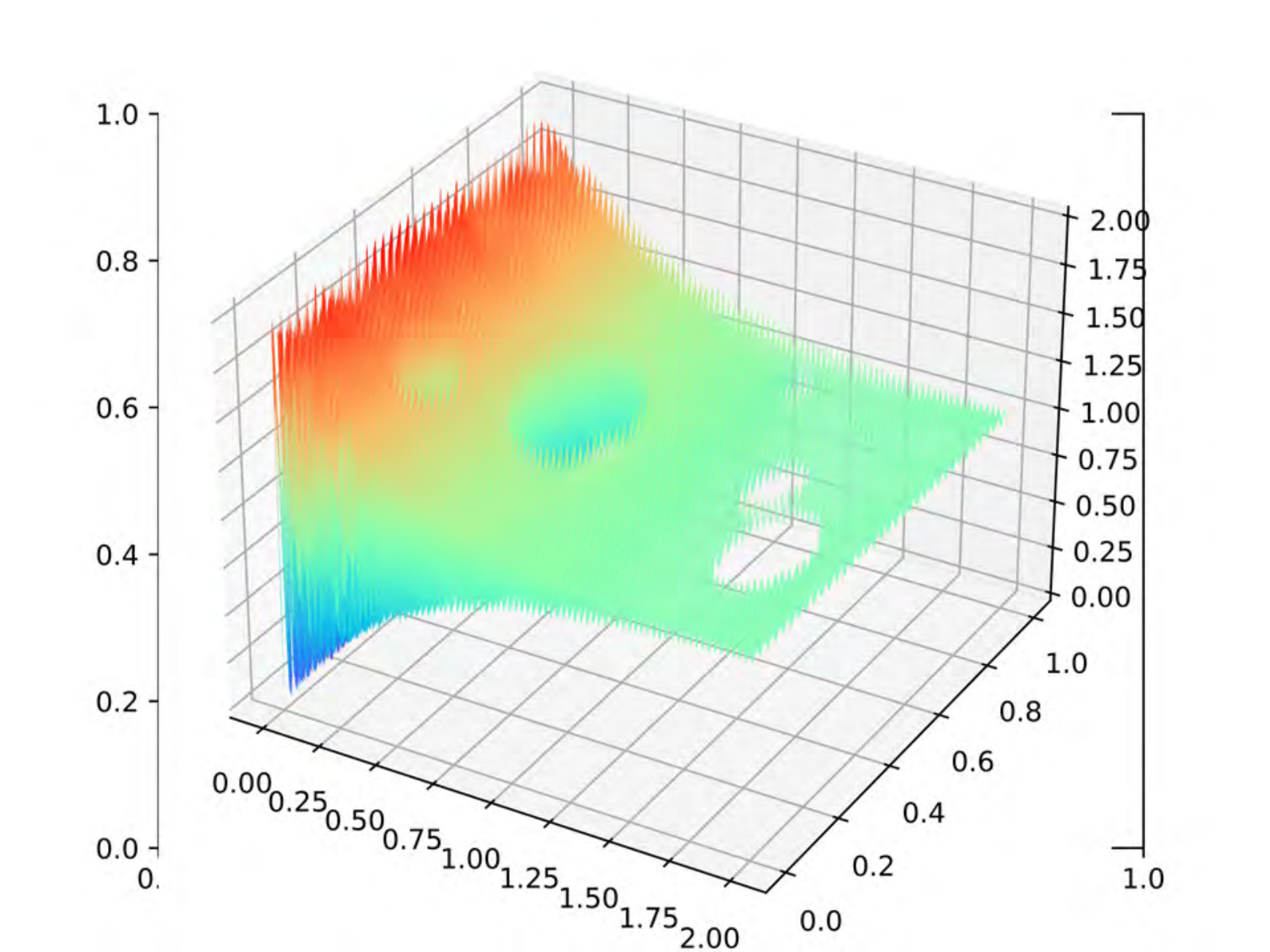}
         \caption{Contour of the first component of velocity of the oscillatory case after 1000 epoch training for scheme Hybrid}
         \label{fig:complexDomain_velocity_x_Contour_Hybrid}
     \end{subfigure}
     \hfill
     \caption{The results of the first component of velocity of the oscillatory case for different schemes}
        \label{fig:complexDomain_velocity_x_Contour}
\end{figure}

\begin{figure}[ptbh]
         \centering
         \includegraphics[width=\textwidth]{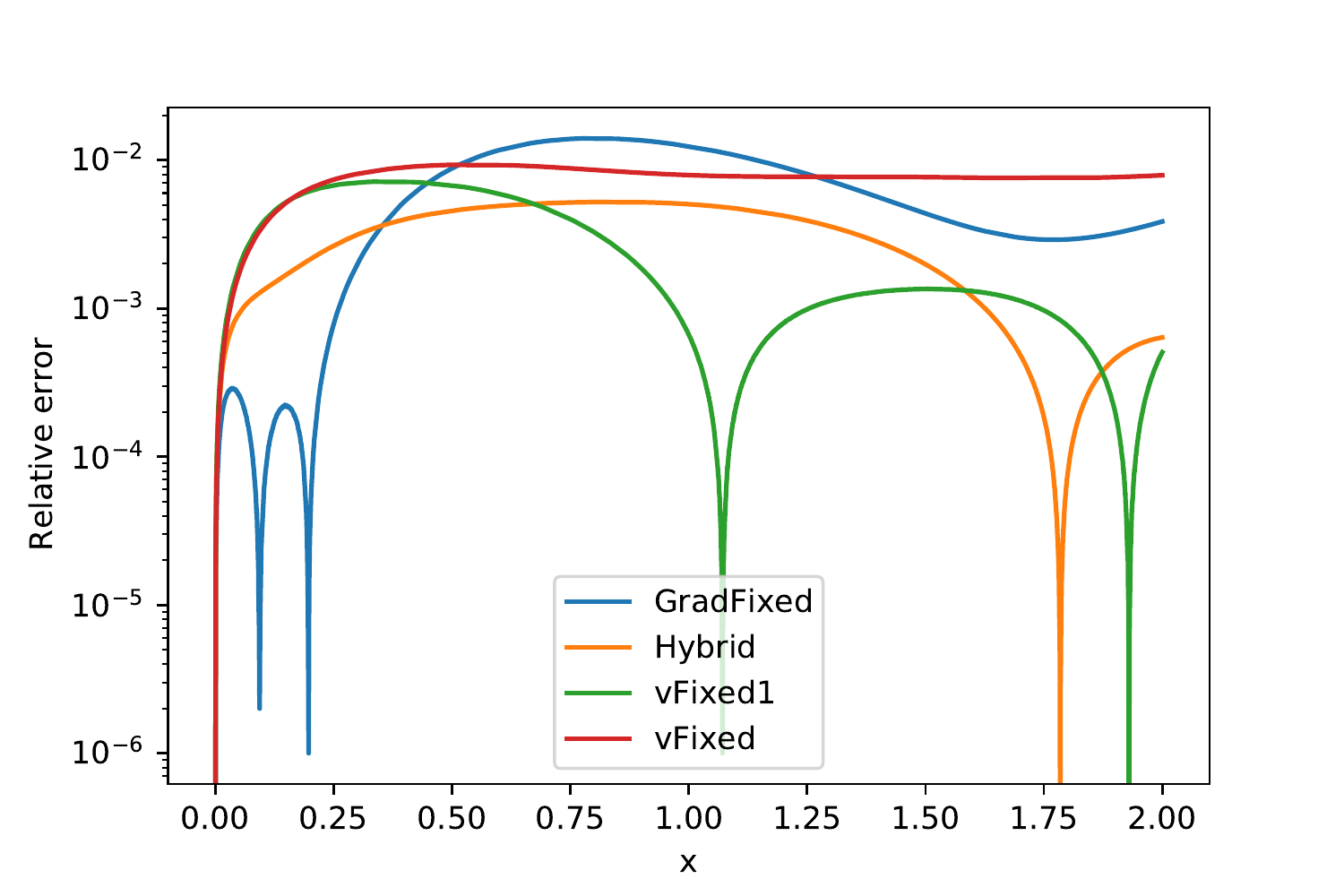}
         \caption{Relative errors alone line \(y = 0.7\) of 4 different schemes for pressure of the complex domain case after 1000 epoch training}
         \label{fig:4scheme_err_p}
\end{figure}

\begin{figure}
     \centering
     \begin{subfigure}[b]{0.48\textwidth}
         \centering
         \includegraphics[width=\textwidth]{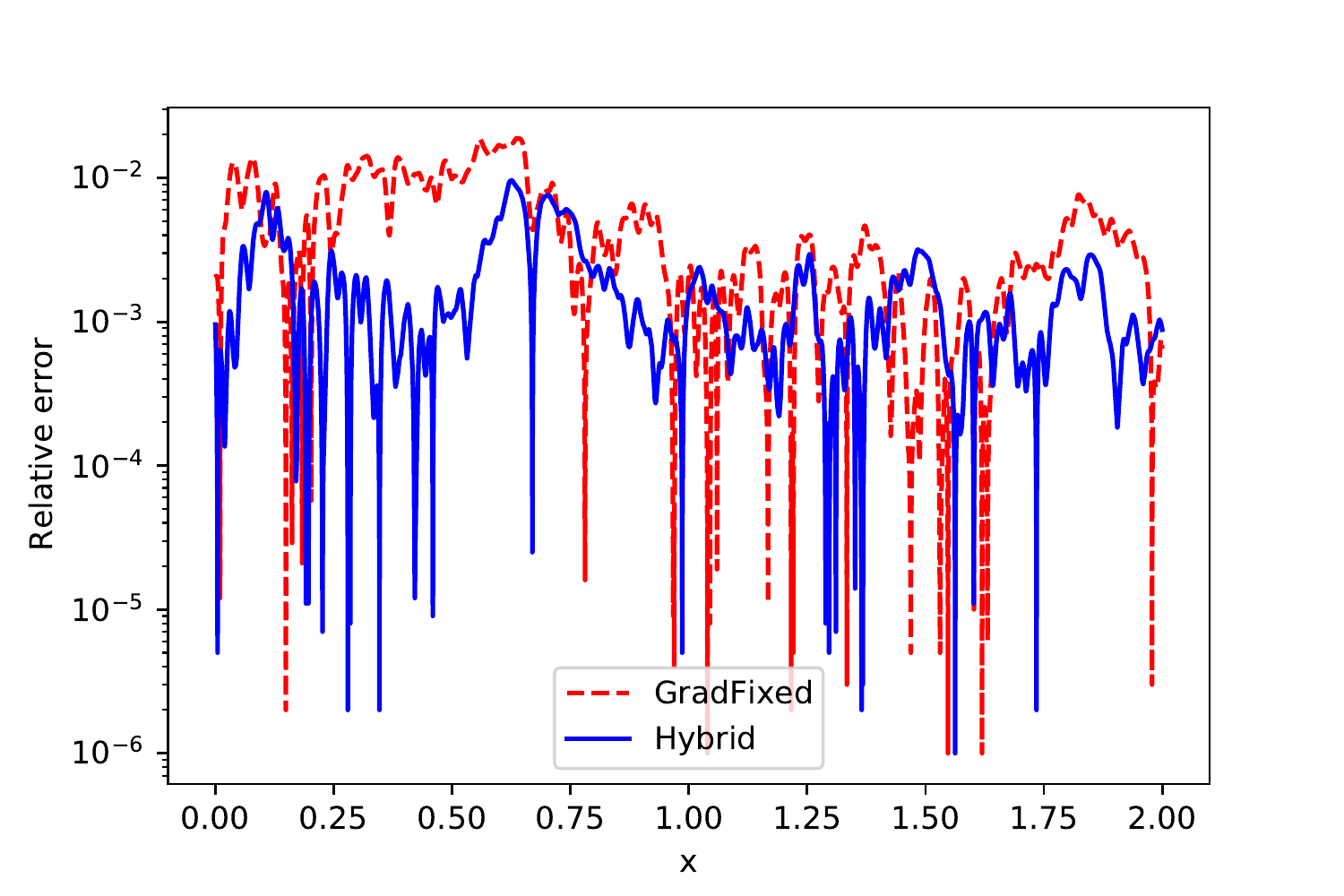}
         \caption{Relative errors alone line \(y = 0.7\) of GradFixed scheme and Hybrid scheme for pressure of the complex domain case after 1000 epoch training}
         \label{fig:4scheme_err_u1}
     \end{subfigure}
     \hfill
     \begin{subfigure}[b]{0.48\textwidth}
         \centering
         \includegraphics[width=\textwidth]{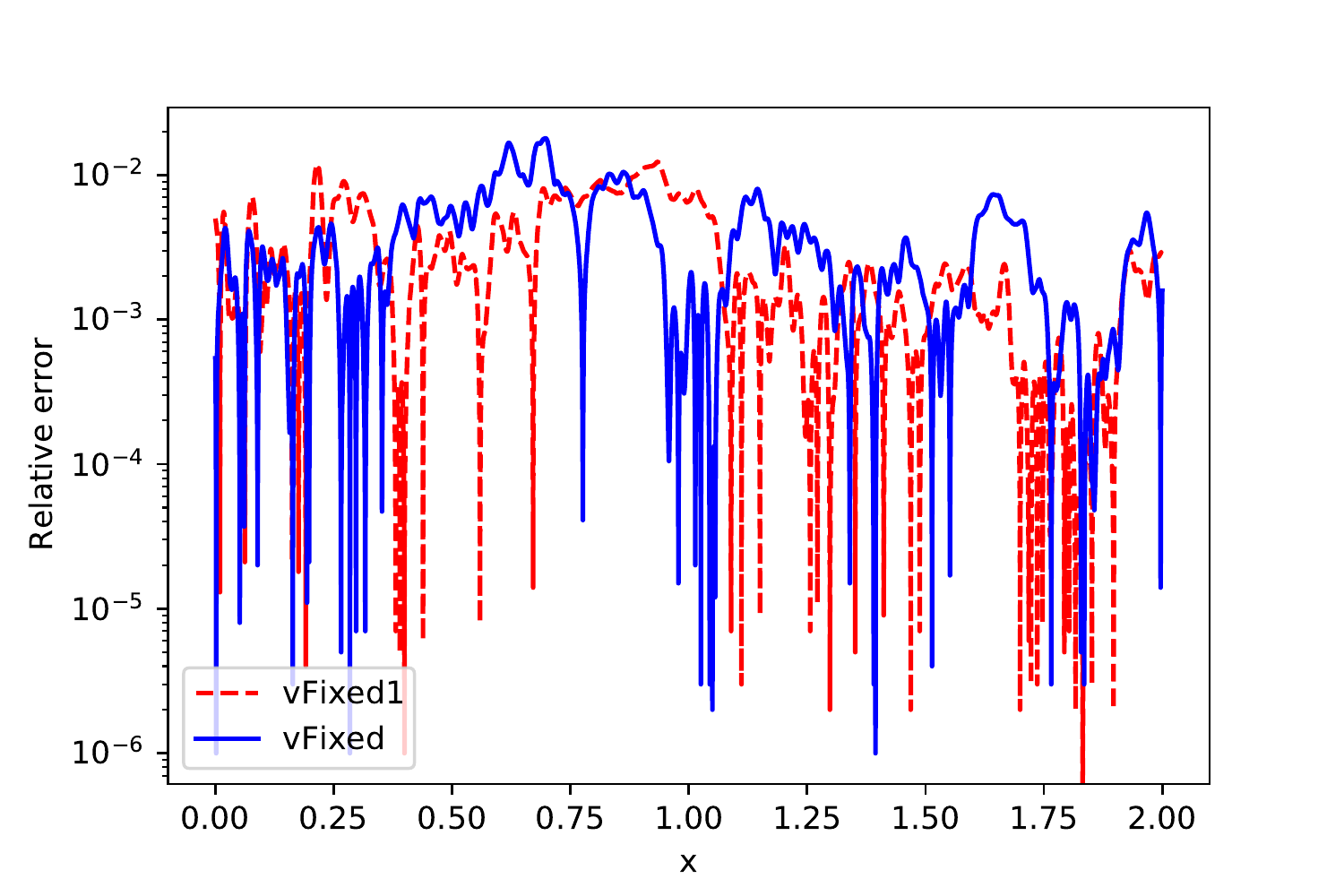}
         \caption{Relative errors along line \(y = 0.7\) of vFixed scheme and vFixed1 for the first component of velocity of the complex domain case after 1000 epoch training}
         \label{fig:4scheme_err_u2}
     \end{subfigure}
     \hfill
     \caption{The relative errors of the complex domain case for different schemes}
        \label{fig:4scheme_err_u}
\end{figure}

In this case, we use the similar settings for the multiscale deep neural networks, adjustments of learning rates, and sampling strategies like what we choose in the oscillatory case with the scale of the simple domain with one hole. The multiscale deep neural networks has 8 scales: \(\{x, 2x, 4x,8x,32x,64x, 128x\}\), whose subnetworks contain 4 hidden layers and 128 hidden neurons in each layer. The frequency we select for this case are the same as what in section \ref{oscillating_simple}. The initial learning rate for this case is \(4e-3\).
Figure \ref{fig:complexDomain_velocity_x_Contour} shows contours of the first component of velocity for different schemes. Figure \ref{fig:4scheme_err_u} and Figure \ref{fig:4scheme_err_p} show the relative errors of these 4 different linearization schemes along line \(y = 0.7\). All schemes we propose converges accurately to the exact solutions.

\section{Conclusion and future work}
\label{sec:future}
In this paper we proposed four linearized learning schemes to solve the stationary highly oscillatory Navier-Stokes flows with multiscale deep neural networks and showed the acceleration of convergence of the schemes are substantial, which demonstrate the capability of the Mscale deep neural networks and the effectiveness of the linearized schemes to solve the nonlinear Navier-Stokes equations. These schemes shed some lights on the practical applications of neural network machine learning algorithms to the nonlinear equations, which are time-consuming using traditional finite element methods. The deep neural network based methods offer an alternative that doesn't require meshes without the need to solve large-scale linear systems, in constrast to traditional numerical methods.

There are much more works to be done for these linearized learning methods, among them the most important is to understand the convergence property of these schemes. The applications of these schemes to other nonlinear PDEs will also be considered. Another challenging work is to consider the time dependent Navier-Stokes equation, which will be explored in the near future.

\medskip
\bibliographystyle{unsrt}

\end{document}